\documentclass[12pt]{article}
\usepackage{amsfonts, amssymb, amsmath, amsthm, latexsym, array}
\usepackage{fullpage}
\usepackage{verbatim}

\usepackage[T2A]{fontenc}
\usepackage[koi8-r]{inputenc}

\input{xy}
\xyoption{all}

\newtheorem{thm}{Theorem}
\newtheorem{lm}{Lemma}
\newtheorem{cl}{Corollary}
\newtheorem{prop}{Proposition}

\theoremstyle{remark}
\newtheorem{ex}{Example}
\newtheorem{rmk}{Remark}

\theoremstyle{definition}
\newtheorem{df}{Definition}

\renewcommand{\iff}{if and only if }

\newcommand{\mymathbb}{\mathrm{{I}\!\!\!1}}

\newcommand{\cled }{\noindent \bf Corollary }

\newcommand{\gt}{\mathfrak}
\newcommand{\cp}{\mathbb C}
\newcommand{\rl}{\mathbb R}
\newcommand{\qv}{\mathbb H}
\newcommand{\qw}{\mathbb H}

\newcommand{\ad}{{\rm ad }}
\newcommand{\SL}{{\rm SL}}
\newcommand{\GL}{{\rm GL}}
\newcommand{\SO}{{\rm SO}}
\newcommand{\Opp}{{\rm O}}
\newcommand{\Upp}{{\rm U}}
\newcommand{\Spp}{{\rm S}}
\newcommand{\SU}{{\rm SU}}

\renewcommand{\U}{{\rm U}}
\newcommand{\Sp}{{\rm Sp}}
\newcommand{\Spin}{{\rm Spin}}

\newcommand{\Lie}{{\rm Lie\,}}
\newcommand{\gr}{{\rm gr\,}}

\begin{document}

\begin{center}{
\Large
\bf
Principal Gelfand pairs}
\par\medskip
{\large Oksana Yakimova}
\end{center}
\vskip2ex

Let $X=G/K$ be a connected Riemannian homogeneous space of
a real Lie group $G$.
We assume that the action $G:X$ of $G$ on $X$ is locally effective, 
i.e., $K$ contains no non-trivial connected central subgroups of $G$.   
Denote by ${\cal D}(X)^G$ the algebra
of $G$-invariant differential operators on $X$ and by ${\cal P}(T^*X)^G$
the algebra of $G$-invariant functions on $T^*X$ polynomial on fibres.
It is well known that ${\cal P}(T^*X)^G$ is
a Poisson algebra, the Poisson bracket being induced by the
commutator in ${\cal D}(X)^G$.

The homogeneous space $X$ is called
{\it commutative} or the pair
$(G, K)$ is called a {\it Gelfand pair} if the following
three equivalent conditions are satisfied:

1) the algebra of $K$-invariant measures on $X$ with
compact support is commutative with respect to the convolution;

2) the algebra ${\cal D}(X)^G$ is commutative;

3) the algebra  ${\cal P}(T^*X)^G$ is commutative with
respect to the Poisson bracket.

The equivalence of the first two conditions was proved by
Thomas \cite{th} and Helgason \cite{h}, independently.
Evidently, condition 3) is a consequence of 2). The
inverse implication is recently
proved by Rybnikov \cite{R}.
Commutative spaces can also be characterised by several other conditions.
For instance,
$X$ is commutative \iff the representation of $G$ in $L^2(X)$ has
a simple spectrum,
see \cite{bg}. 

Symmetric Riemannian homogeneous spaces introduced by \'Elie Cartan
are commutative.
The theory of symmetric spaces is well developed.
Works of \'Elie Cartan and Helgason
describe their structure and also deal with harmonic
analysis on such manifolds. One can hope that some day 
commutative spaces will be as thoroughly studied as symmetric spaces. 

We denote Lie algebras of Lie groups by
corresponding small gothic letters, for example, $\gt g=\Lie G$
and $\gt k=\Lie K$. 

\begin{df} A real or complex linear Lie group
with finitely many connected components is said to be
{\it reductive} if
it is completely reducible.
\end{df}

If $X$ is a commutative homogeneous space of $G$, then, up
to the local isomorphism, $G$ has a factorisation
$G=N\leftthreetimes L$, where $N$ is a nilpotent radical of $G$,
$K\subset L$, $L$ and $K$ have the
same invariants in $\rl[\gt n]$ 
and $[[\gt n, \gt n], \gt n]=0$, see \cite{Vin}. Without loss of
generality we may assume that the center of $L$ is compact
and the commutator group $L'$ of $L$ is a real form of some complex 
semisimple group.
Hence, $L$ is a reductive group.
 For a reductive Lie group $F$ acting
on a linear space $V$, we denote by $F_*(V)$ a generic
stabiliser for this action. The
 stabiliser of a point $y\in V$ is denoted by $F_y$.
Set 
$\gt m:=\gt l/\gt k$.
In section 1 we prove the following criterion of commutativity.

\vskip0.2ex 
\noindent 
{\bf Theorem 0.1.} {\it
$X=(N\leftthreetimes L)/K$ is commutative if and only if 
all of the following three conditions hold:

i) $\rl[\gt n]^L=\rl[\gt n]^K$;

ii) for any point $\gamma\in\gt n^*$ the
homogeneous space $L_\gamma/K_\gamma$ is commutative;

iii) for any point $\beta\in\gt m^*$ the homogeneous space
$(N\leftthreetimes K_\beta)/K_\beta$ is commutative.
}
\vskip1ex

Let $F$ be a complex reductive Lie group
and $H\subset F$ a reductive subgroup.

\begin{df}
An  affine complex $F$-variety $X$ is called {\it spherical} if
a Borel subgroup $B(F)\subset F$ has an open orbit in $X$.
If $X$ is a linear space
and a spherical $F$-variety then it is called
a spherical representation of $F$. If a homogeneous space
$F/H$ is spherical, then
the pair $(F, H)$ and the subgroup $H$ are also called spherical.
\end{df}

Let $G$ be a real
form of a complex reductive group
$G(\cp)$. Suppose $K\subset G$ is a compact subgroup.
We call the real homogeneous space $G/K$, the subgroup $K$ and
the pair $(G, K)$ spherical if
the complexification $X(\cp)=G(\cp)/K(\cp)$ is a spherical
$G(\cp)$-variety.

The commutative homogeneous spaces of reductive Lie groups are
just the real forms of the spherical
affine homogeneous spaces, see, for example, \cite{Vin}. The theory of
spherical homogeneous spaces is
well developed, in particular, they are classified in
 \cite{kr},
\cite{br} and \cite{m}. Note that \cite{br} deals only with
so-called {\it principal} homogeneous spaces (see Definition
5 in the third part of the present article). In \cite{m} one
type of non-principal
spherical homogeneous spaces is described. 
The final classification is obtained in \cite{y2}.
The real forms of homogeneous spherical spaces, i.e.,
commutative homogeneous spaces
are explicitly described in
\cite{y2}.

Let $X=(N\leftthreetimes L)/K$ be a commutative homogeneous space.
Denote by $P$ the ineffective kernel of the action $L:\gt n$.
Evidently,
$P$ is a normal subgroup of $L$ and $G$.
Due to (i) we have $L/P\subset\Opp(\gt n)$.
Hence the group $L_\gamma$ is reductive for any
$\gamma\in\gt n^*$, so in (ii) we can replace the
word ``commutative'' by ``spherical'' and look up the given
homogeneous space in the list.

Because the orbits of the compact group $K$ in $\gt n$ are separated
by polynomial invariants, $L$ and $K$ have the same
invariants in $\rl [\gt n]$ if and only if they have the same
orbits. In other words, condition (i) means that there is a
factorisation $L=L_*(\gt n)K$ or equivalently
$L/P$ is a product of $L_*(\gt n)/P$ and $K/(K\cap P)$. 
All nontrivial factorisations of
compact groups into products of two subgroups are
classified by Onishchik \cite{On}. The classification
results are also nicely presented in  \cite[Chapter 4]{On2}.

Commutativity is a local property, i.e., it depends only on
the pair of algebras $(\gt g, \gt k)$, see \cite{Vin}. We
assume that both $G$ and $K$ are connected, $N$ is simply
connected, $L=Z(L)\times L_1\times...\times L_m$, where
$Z(L)$ stands for the connected centre of $L$ and $L_i$ are
connected non-commutative simple groups. We also assume that
$L_i$ are real forms of simply connected complex simple
groups and the action of a factor group $Z(L)/(Z(L)\cap
P^0)$ on $\gt n$ is effective. It can happen, that for a
given pair $(\gt g, \gt k)$, there is no effective pair $(G,
K)$ satisfying these assumptions, so we admit not only
effective actions $G:(G/K)$ but locally effective as well.
In tables and theorems we write $U_n$ instead of $U_1\times\SU_n$ and
sometimes $\SO_n$ instead of $\Spin_n$.

Suppose $X=(N\leftthreetimes L)/K$ is commutative.
Let $\gt z_0\subset[\gt n,\gt n]$ be an 
$L$-invariant subspace 
and $Z_0\subset N$ a corresponding connected subgroup. Then the
homogeneous space $X/Z_0=((N/Z_0)\leftthreetimes L)/K$ is also commutative,
see \cite{Vin}.
The passage from $X$ to $X/Z_0$ is called a {\it central reduction}.
A commutative homogeneous space is said to be {\it maximal},
if it can not be obtained by a
non-trivial central reduction from a larger one.
We consider only maximal commutative homogeneous spaces. 
All commutative homogeneous spaces could be obtained 
as their central reductions.  

\begin{df} A homogeneous space $G/K$ is called {\it indecomposable}
if it cannot be represented as a product
$G_1/K_1\times G_2/K_2$, where $G=G_1\times G_2$,
$K=K_1\times K_2$ and $K_i\subset G_i$.
\end{df}

Obviously, $G_1/K_1\times G_2/K_2$ is commutative
\iff both spaces $G_1/K_1$ and $G_2/K_2$ are commutative.
Hence, for the classification of commutative homogeneous spaces
it is enough to describe indecomposable ones.

Denote by $H_n$ the $2n{+}1$ dimensional Heisenberg group, 
i.e., $\gt h_n=\Lie H_n\cong\cp^n\oplus\rl$.  
Simply connected commutative groups are denoted by 
$\rl^n$ or $\cp^n$. 
The simplest and most important results are obtained 
for simple $L$.  

\vskip0.2ex
\noindent
{\bf Theorem 0.2.} {\it
Suppose $X=(N\leftthreetimes L)/K$ is an indecomposable
commutative space, where $L$ is simple, $L\ne K$ and 
$\gt n\ne 0$. Then $X$ is one of the following eight spaces.
}
$$
\begin{array}{lll}
(H_{2n}\leftthreetimes\SU_{2n})/\Sp_n; \qquad &
   (\rl^7\leftthreetimes\SO_7)/G_2; \qquad & 
     ((\rl^8\times\rl^2)\leftthreetimes\SO_8)/\Spin_7; \\
(\cp^{2n}\leftthreetimes\SU_{2n})/\Sp_n; \qquad &
     (\rl^8\leftthreetimes\Spin_7)/\Spin_6; \qquad & 
           (\rl^8\leftthreetimes\SO_8)/\Spin_7;\\
(\rl^{2n}\leftthreetimes\SO_{2n})/\SU_n; \qquad & &
         (\rl^8\leftthreetimes\SO_8)/(\Sp_2\times\SU_2).\\ 
\end{array}
$$
\vskip02.ex

First we describe commutative homogeneous spaces satisfying 
condition
\vskip0.5ex
\hbox to \textwidth{\qquad \quad $({\bf \ast})$\hfil $L\ne
K$ and
the action $L:\gt n$ is locally effective, i.e.,
$P$ is finite.
\hfil }

\vskip0.2ex
\noindent
{\bf Theorem 0.3.} {\it
Suppose a commutative homogeneous space $X=G/K$ satisfies 
condition~$({\bf \ast})$.
Then any non-commutative normal subgroup of $K$ different  from
$\SU_2$ 
is contained in some simple direct factor of $L$.
}
\vskip02.ex

In the present work we impose on $X$ two technical conditions of
``principality'' and ``$\Sp_1$-saturation'', see sections 3 and 5
for precise definitions
and explanations. The first condition concerns the embedding of
the connected centre of $K$ into $L$ and the action of the connected
centre of $L$ on $\gt n$. 
The second condition describes the behaviour of normal subgroups of
$K$ and $L$ isomorphic to $\Sp_1$. 
Example~\ref{derevo} in the beginning 
of Section~5 shows, that the classification in the 
general case requires 
a lot of tedious calculations. 

\vskip0.3ex
\noindent
{\bf Theorem 0.4.} {\it
Let $X=(N\leftthreetimes L)/K$ be a maximal
principal indecomposable
commutative homogeneous space satisfying
condition $({\bf \ast})$. Suppose there is a 
simple normal subgroup $L_1\lhd L$, such that 
$L_1\ne\SU_2$ and $L_1\not\subset K$. 
Then either $L$ is simple (and $X$ is listed in Theorem 0.2) or 
$X$ is one of the following two spaces:
$(H_{2n}\leftthreetimes\U_{2n})/(\Sp_n\cdot\U_1)$, 
$((\rl^8\times\rl^2)\leftthreetimes(\SO_8\times\SO_2))/(\Spin_7\times\SO_2)$.}
\vskip0.1ex

These additional spaces are essentially the same as 
the spaces from Theorem 0.2. All of them are listed in Table 2b (section 2). 

\vskip0.1ex

A  commutative homogeneous space $(N\leftthreetimes L)/K$ 
is said to be of 
{\it Heisenberg type} if $L=K$. 
Recently commutative homogeneous spaces of Heisenberg type 
were intensively studied by
several people, see, e.g., \cite{be-ra}, \cite{la}, 
\cite{nish}, \cite{Vin}, \cite{Vin2}. 

The following theorem is the main result of this article.

\vskip0.3ex
\noindent
{\bf Theorem 0.5.} {\it Any indecomposable 
maximal principal $\Sp_1$-saturated commutative homogeneous 
space belongs to the
one of the following 
four classes:

\noindent
1) affine spherical homogeneous spaces of reductive real Lie groups;

\noindent
2) homogeneous spaces listed in Table 2b;

\noindent
3) homogeneous spaces $((\rl^n\leftthreetimes\SO_n)\times\SO_n)/\SO_n$,
$((H_n\leftthreetimes\U_n)\times\SU_n)/\Upp_n$,
where the normal subgroups $\SO_n$ and $\SU_n$ of $K$
are diagonally embedded into $\SO_n\times\SO_n$ and
$\SU_n\times\SU_n$, respectively;

\noindent
4) commutative  homogeneous spaces of Heisenberg type.
}\vskip0.3ex


Commutative homogeneous spaces of Heisenberg type
are considered in the sixth section. 
In this case ${\cal D}(G/K)^G\cong U(\gt n)^K$, 
where $U(\gt n)$ is the universal enveloping algebra 
of $\gt n$. 
If $\gt n$ is commutative then obviously $G/K$ is also commutative and it is
called a commutative space of {\it Euclidian type}. We assume that
$\gt n$ is not commutative.

Consider a homogeneous space $(N\leftthreetimes K)/K$. 
Suppose $\gt n$ is at most two-step nilpotent and 
$[\gt n, \gt n]$ is a trivial $K$-module.
We can decompose $\gt n$ into an $K$-invariant sum
$\gt n=(\gt w\oplus\gt z)\oplus V$, where
$V$ is an abelian ideal and $[\gt w,\gt w]=\gt z$. 
Any point $\alpha\in\gt z^*$ determines a 
skew-symmetric form $\hat\alpha$ on $\gt w$, 
namely $\hat\alpha(\xi,\eta)=\alpha([\xi,\eta])$ 
for $\xi,\eta\in\gt w$. The form $\hat\alpha$ 
is non-degenerate for a generic $\alpha$. 
The complexification $\gt w(\cp)$ is an orthogonal 
and a symplectic representation of $K(\cp)$ at the same time. 
Hence it is reducible $\gt w(\cp)=W\oplus W^*$.  
According to \cite{be-ra} and \cite{Acta}, 
$(N\leftthreetimes K)/K$ is commutative \iff
$W$ is a spherical representation of the complexification
$K_*(V)(\cp)$ of $K_*(V)$.
In the simplest situation when $V=0$ the statement means that $W$
is a spherical representation of $K(\cp)$.

Spherical representations of reductive Lie groups
are classified
by Kac \cite{Kac} (irreducible representations), 
Brion \cite{br2}, finally by  Benson and Ratcliff \cite{be-ra2} 
and Leahy \cite{L}, independently.
Historical comments and the classification result can be found
in \cite{kn}.
The list of commutative homogeneous spaces 
$(N\leftthreetimes K)/K$, 
where $N$ is a product of several Heisenberg groups,
is given in \cite{be-ra}. That article also claims to classify 
all commutative homogeneous spaces 
$(N\leftthreetimes K)/K$ with  
$\gt n=(\gt w\oplus\gt z)\oplus V$, where 
$V$ is an abelian ideal and $\gt w\oplus\gt z$
is a direct sum of several Heisenberg algebras.  
The authors of \cite{be-ra} erroneously assume
that if for a subgroup $H$ with
$\Lie H=\gt w\oplus\gt z$ the homogeneous space
$(H\leftthreetimes K)/K$ is commutative, then  
$(N\leftthreetimes K)/K$ is also commutative. 
The simplest counterexample is 
$((\cp^2\times H_2)\leftthreetimes\SU_2)/\SU_2$. 
This space is not commutative according to \cite{Vin}, but
$(H_2\leftthreetimes\SU_2)/\SU_2$ is commutative.

Commutative homogeneous spaces of Heisenberg type with 
an irreducible action
$K:(\gt n/\gt n')$ are classified in \cite{Vin} and \cite{Vin2}.
Generally, $\gt n$ is a sum of a commutative ideal and algebras listed in
\cite[Table 3]{Vin} and \cite{Vin2}.
But the sum cannot be arbitrary.

The classification of {\it saturated} commutative
homogeneous spaces of Heisenberg type
was announced in \cite{Acta}. The condition of saturation
is a little bit stronger than both conditions of
``principality'' and ``$\Sp_1$-saturation''. We
present the result of \cite{Acta} in Table 4 and give
a proof.

\par\smallskip

\noindent

We will frequently use the following results of \cite{Vin}.
\begin{prop}{\rm \cite[Corollaries to Proposition 10]{Vin}}
Let $G/K$ be commutative. Then

1) for any normal subgroup $N\subset G$ the homogeneous space $G/NK=
(G/N)/(K/(N\cap K))$ is commutative;

2) for any compact subgroup $F\subset G$ containing $K$
the homogeneous space $G/F$ is commutative;

3) for any subgroup $F\subset G$
containing $K$  the homogeneous space $F/K$ is  commutative.
\end{prop}

\vskip0.1ex
\noindent
We fix some additional notation.

$G'$ is the commutator group of $G$;

$G(\cp)$ is the complexification of a real group $G$;

$B(G(\cp))$ and $U(G(\cp))\subset B(G(\cp))$ are a Borel and a maximal
unipotent subgroups of a reductive group $G(\cp)$;

$X/\!\!/G$ stands for the categorical quotient of
an affine algebraic variety $X$ by the action of a reductive 
group $G$.

\par\medskip

{\bf Acknowledgments.} 
I am grateful to 
E.B.\,Vinberg and D.I.~Panyushev for helpful discussions
and permanent attention to this work. 
Thanks are also due to A.L.~Onishchik and D.N.~Akhiezer for 
useful comments on the earlier version. 
This paper was finished during my stay at the Max-Planck-Institut
f\"ur Mathematik (Bonn). I would like to thank 
the Institut for hospitality and wonderful
working conditions. This research was supported in part by 
CRDF grant no. RM1-2543-MO-03. 

\par\bigskip

\begin{center}
{\large \bf 1. Commutativity criterion}
\end{center}

Let $U(\gt g)$ stand for the universal
enveloping algebra of $\gt g$.
There is a natural filtration:
$$
  U_0(\gt g)\subset U_1(\gt g)\subset \ldots \subset
                U_m(\gt g)\subset\ldots \ ,
$$
where $U_m(\gt g)\subset U(\gt g)$
consists of all elements of order at most $m$.

The Poisson bracket on the symmetric algebra
 $S(\gt g)=\gr U(\gt g)$ is determined by the formula
$$
 \{a+U_{n-1}(\gt g), b+U_{m-1}(\gt g)\}=[a, b]+U_{n+m-2}(\gt g)
   \quad \forall
a\in U_n(\gt g), b\in U_m(\gt g).
$$
Let $X=G/K$ be a Riemannian homogeneous space. 
It is well known, see, for example, \cite{Vin}, that
there is an isomorphism of the associated graded algebras:
$$
\gr U(\gt g)^K/(U(\gt g)\gt k)^K=\gr{\cal D}(X)^G={\cal P}(T^*X)^G=
S(\gt g/\gt k)^K.
$$
The space $(U(\gt g)\gt k)^K$ is an ideal of $U(\gt g)^K$,
also $(S(\gt g)\gt k)^K$ is a Poisson ideal of $S(\gt g)^K$.
The well defined Poisson bracket on the
factor $S(\gt g)^K/(S(\gt g)\gt k)^K\cong S(\gt g/\gt k)^K$
coincides up to a sign with the Poisson bracket on
${\cal P}(T^*X)^G$. In particular, $X$ is commutative \iff
the Poisson algebra $S(\gt g/\gt k)^K$ is commutative.

Suppose $X=(N\leftthreetimes L)/K$ 
is commutative. Then, as proved in \cite{Vin},
the following condition holds

i) $\rl[\gt n]^L=\rl[\gt n]^K$.

\noindent
The orbits of a compact group are separated by polynomial invariants. Hence
(i) is fulfilled \iff $L$ and $K$ have the same orbits on $\gt n$.
There is a $K$-invariant positive-definite symmetric bilinear form on $\gt n$
which is automatically $L$-invariant. In particular, vector spaces $\gt n$
and $\gt n^*$ are
isomorphic as $L$-modules. Therefore, $\ad^*(\gt k)\gamma=\ad^*(\gt l)\gamma$
for each point $\gamma\in\gt n^*$ and hence $\gt l=\gt k+\gt l_\gamma$.
Moreover, the natural restriction
$$
\tau:\enskip \gt l^*\longrightarrow\gt l_{\gamma}^*
$$
(which is also a homomorphism of $L_\gamma$-modules) determines an
isomorphism of
$K_{\gamma}$-modules $(\gt l/\gt k)^*$ and  $(\gt
l_{\gamma}/\gt k_{\gamma})^*$.

Recall that $\gt g=\gt l+\gt n$, where $\gt n$ is a
nilpotent ideal and $\gt l$
is a reductive subalgebra.
Let $\check{\gt n}$ and $\check{\gt l}$ be commutative Lie algebras
of dimensions $\dim\gt n$ and $\dim\gt l$,
respectively.
We determine new Lie algebras
$\check{\gt g}_1=\gt l +\check{\gt n}$ and $\check{\gt
g}_2=\check{\gt l}\oplus\gt n$, where
$\check{\gt l}$, $\check{\gt n}$ are commutative ideals
and $\check{\gt n}\cong\gt n$
as an $\gt l$-modules.

Denote by $\{\phantom{.}, \phantom{.}\}_{\gt l}$
and $\{\phantom{.}, \phantom{.}\}_{\gt n}$
the Poisson brackets on
$S(\check{\gt g}_1)$ and $S(\check{\gt g}_2)$.
There is a $K$-invariant bi-grading
$S(\gt g)=\bigoplus
S^{n,l}(\gt g)$, where $S^{n,l}(\gt g)=S^n(\gt n)S^l(\gt l)$.
We identify elements of $S(\gt g)$ with the
corresponding elements of $S(\check{\gt g}_1)$
and $S(\check{\gt g}_2)$.

\begin{lm}\label{vs}
We have
$$ \{\xi, \eta\}=\{\xi, \eta\}_{\gt n}+\{\xi, \eta\}_{\gt l},
\enskip
\mbox{ with } \{\xi, \eta\}_{\gt n}\in S^{n+n'-1, l+l'}(\gt g),
\{\xi, \eta\}_{\gt l}\in S^{n+n', l+l'-1}(\gt g).
$$
for any bi-homogeneous elements
$\xi\in S^{n,l}(\gt g), \eta\in S^{n', l'}(\gt g)$.
In other words,
the Poisson bracket on $S(\gt g)$ is a direct sum of the
brackets
$\{\phantom{.}, \phantom{.}\}_{\gt n}$ and
$\{\phantom{.}, \phantom{.}\}_{\gt l}$.
\end{lm}

\begin{proof}
The Poisson bracket of bi-homogeneous
elements
$\xi=\xi_1...\xi_n, \eta=\eta_1...\eta_m\in S(\gt g)$
is given by the formula
\begin{equation}\label{pbr}
\{\xi, \eta\}=\sum\limits_{i,j}[\xi_i,\eta_j]
\xi_1...\widehat{\xi_i}...\xi_n\eta_1...\widehat{\eta_j}...\eta_m.
\end{equation}
This expression for $\{\xi, \eta\}$ contains summands of three
different types, depending on whether $\xi_i$
and $\eta_j$ are elements of $\gt l$ or $\gt n$.
Because $[\gt l, \gt n]\subset\gt n$ and $\gt l$,
$\gt n$ are subalgebras, if
$\xi_i, \eta_j\in\gt n$, then
$[\xi_i, \eta_j]\in S^{n+n'-1, l+l'}(\gt g)$,
otherwise $[\xi_i, \eta_j]\in S^{n+n', l+l'-1}(\gt g)$.
\end{proof}

We suppose that $\check{\gt k}$ is embedded into
$\check{\gt l}$ as a commutative
subalgebra of dimension $\dim\gt k$.
We also denote by $\{\phantom{.}, \phantom{.}\}_{\gt n}$ and
$\{\phantom{.}, \phantom{.}\}_{\gt l}$ the corresponding
Poisson brackets on the Poisson factors
$S(\check{\gt g}_2/\check{\gt k})^K=
S(\check{\gt g}_2)^K/(S(\check{\gt g}_2)\check{\gt k})^K$
and $S(\check{\gt g}_1/\gt k)^K=
S(\check{\gt g}_1)^K/(S(\check{\gt g}_1)\gt k)^K$,
where the actions $K:\check{\gt g}_i$ are the same as
$K:\gt g$. We have
$ \{a, b\}_{\gt l}\in S^{n+n', l+l'-1}(\gt g/\gt k)$
and
$\{a, b\}_{\gt n}\in S^{n+n'-1, l+l'}(\gt g/\gt k)$
for any $a\in S^{n,l}(\gt g/\gt k), b\in S^{n', l'}(\gt g/\gt k)$
($a,b\in S(\gt g/\gt k)^K$).

\begin{lm}\label{kr1} The Poisson bracket on
$S(\gt g/\gt k)^K$ is of the form
$\{\phantom{.}, \phantom{.}\}=\{\phantom{.},
\phantom{.}\}_{\gt n}+\{\phantom{.}, \phantom{.}\}_{\gt l}$.
\end{lm}

\begin{proof} This is  a straightforward  consequence of Lemma \ref{vs}.
\end{proof}

\begin{cl} Let $G/K$ be a commutative homogeneous space and $\check N$ a
simply connected commutative Lie
group with a Lie algebra $\check{\gt n}$. Set
$\check G:=\check N\leftthreetimes L$.
Then $\check G/K$ is also commutative.
\end{cl}

\begin{thm} The homogeneous space
$X=(N\leftthreetimes L)/K$ is  commutative if and only if 
all of the following three
conditions hold:

i) $\rl[\gt n]^L=\rl[\gt n]^K$;

ii) for any point $\gamma\in\gt n^*$ the
homogeneous space $L_\gamma/K_\gamma$ is commutative;

iii) for any point $\beta\in(\gt l/\gt k)^*$ the homogeneous space
$(N\leftthreetimes K_\beta)/K_\beta$ is commutative.
\end{thm}

\begin{rmk} The statement of the theorem remains true if we replace
 arbitrary points by generic points
in conditions (ii) and (iii).
\end{rmk}

\begin{proof}
As was already mentioned, Vinberg proved in \cite{Vin} that the
condition  (i) holds for any commutative space. So let us
assume that it is fulfilled.

Let $\gamma$ be a point in $\gt n^*$.
Recall that the $K_\gamma$-modules $\gt l/\gt k$ and $\gt
l_\gamma/ \gt k_\gamma$ are isomorphic. Hence, $S(\gt l/\gt k)$ is
isomorphic to $S(\gt l_\gamma/\gt k_\gamma)$ as a graded
associative algebra and also as a $K_\gamma$-module.

Consider the homomorphism
$$
\varphi_{\gamma}:\enskip S(\gt g/\gt k)\longrightarrow S(\gt g/\gt
k)/(\xi-\gamma(\xi):\enskip\xi\in\gt n)=S(\gt l/\gt k)=S(\gt l_{\gamma}/
\gt k_{\gamma}).
$$
Evidently, $\varphi_\gamma(S(\gt g/\gt k)^K)\subset
S(\gt l_\gamma/\gt k_\gamma)^{K_{\gamma}}$.

Let $\xi\in\gt l_\gamma$, $\eta\in\gt n$.
Then $\gamma(\{\xi, \eta\})=\gamma([\xi,
\eta])=-[\ad^*(\xi)\gamma](\eta)=0=\{\xi, \gamma(\eta)\}$.

It can easily be deduced from the above statement and from
the formula (\ref{pbr}),  that for arbitrary bi-homogeneous
elements $a,b\in S(\gt g/\gt k)^K$, which can be regarded as elements of
$S((\gt l_\gamma\oplus\gt n)/\gt k_\gamma)$, we have
$$
\varphi_\gamma(\{a, b\}_{\gt l})=\{\varphi_\gamma(a),
\varphi_\gamma(b)\},
$$
where the second bracket is the Poisson bracket
on $S(\gt l_\gamma/\gt k_\gamma)^{K_\gamma}$.
In other words, $\varphi_\gamma$ is a homomorphism of the
Poisson algebras $S(\check{\gt g}_1/\gt k)^K$ and
$S(\gt l_\gamma/\gt k_\gamma)^{K_\gamma}$.

Recall that $\gt m=(\gt l/\gt k)$. 
We repeat the procedure for the point $\beta\in\gt m^*$.
Consider the homomorphism
$$
\varphi_{\beta}:\enskip S(\gt g/\gt k)\longrightarrow S(\gt g/\gt
k)/(\xi-\beta(\xi):\enskip\xi\in\gt m)=S(\gt n).
$$
Clearly, $\varphi_{\beta}(S(\gt g/\gt k)^K)\subset
S(\gt n)^{K_{\beta}}$. Note that $\varphi_\beta$ is
a homomorphism of Poisson algebras
$S(\check{\gt g_2}/\check{\gt k})^K$ and $S(\gt n)^{K_\beta}$.
For arbitrary bi-homogeneous
elements $a,b\in S(\gt g/\gt k)^K$ we have
$$
\varphi_\beta(\{a, b\}_{\gt n})=\{\varphi_\beta(a),
\varphi_\beta(b)\}.
$$
Here the second bracket is a Poisson bracket on $S(\gt n)^{K_\beta}$.

Now we show that homomorphisms $\varphi_\gamma$ and $\varphi_\beta$ are
surjective.
We have $S(\gt g)=\mathbb R[\gt g^*]$, $S(\gt
g/\gt k)^K=\mathbb R[(\gt g/\gt k)^*]^K=\mathbb R[(\gt g/\gt  k)^*/\!\!/K]$ and
$S(\gt l_\gamma/\gt k_\gamma)^{K_\gamma}=\mathbb R[\gt
m^*/\!\!/K_\gamma]$, $S(\gt n)^{K_\beta}=\rl [\gt
n^*/\!\!/K_\beta]$. Note that
$$
 \begin{array}{l}
  \gt m^*/\!\!/K_{\gamma}\cong(K\gamma\times\gt
                           m^*)/\!\!/K\subset (\gt g/\gt  k)^*/\!\!/K;\\
\gt n^*/\!\!/K_{\beta}\cong(\gt n^*\times
                          K\beta)/\!\!/K\subset (\gt g/\gt  k)^*/\!\!/K.\\
\end{array}
$$
Moreover,  $K\gamma$ and  $K\beta$ are closed in $\gt
n^*$ and $\gt m^*$, respectively. Hence the subsets
$(K\gamma\oplus\gt m^*)/\!\!/K$ and
 $(\gt n^*\oplus K\beta)/\!\!/K$ are closed in $(\gt g/\gt  k)^*/\!\!/K$.
Thus, the restrictions $\rl [(\gt g/\gt  k)^*]^K\rightarrow
\mathbb R[K\gamma\oplus\gt m^*]^K$ and  $\rl[(\gt g/\gt  k)^*]^K\rightarrow
\mathbb R[\gt n^*\oplus K\beta]^K$
are surjective. It is therefore proved that $\varphi_\gamma$ and
$\varphi_\beta$
are surjective.

\par\smallskip

\noindent
($\Longleftarrow$) Suppose conditions (ii) and (iii) are satisfied.
Clearly, $X$ is commutative \iff both Poisson brackets
$\{\phantom{.}, \phantom{.}\}_{\gt n}$ and
$\{\phantom{.}, \phantom{.}\}_{\gt l}$ equal zero on $S(\gt g/\gt
k)^K$. If $\{a, b\}_{\gt l}\ne
0$ for some elements $a,b\in S(\gt g/\gt k)^K$ then there is a
(generic) point $\gamma\in\gt n^*$ such that $\varphi_\gamma(\{a,
b\}_{\gt l})\ne 0$. But $\varphi_\gamma(\{a, b\}_{\gt
l})=\{\varphi_\gamma(a), \varphi_\gamma(b)\}=0$. Analogously, if
$\{a, b\}_{\gt n}\ne
0$ for some elements $a,b\in S(\gt g/\gt k)^K$ then there is a
(generic) point $\beta\in\gt m^*$ such that $\varphi_\beta(\{a,
b\}_{\gt l})\ne 0$. But $\varphi_\beta(\{a, b\}_{\gt
l})=\{\varphi_\beta(a), \varphi_\beta(b)\}=0$.

\par\smallskip

\noindent
($\Longrightarrow$) Suppose $X$ is commutative. Then both Poisson brackets
$\{\phantom{.}, \phantom{.}\}_{\gt
n}$ and  $\{\phantom{.}, \phantom{.}\}_{\gt l}$ vanish on $S(\gt g/\gt k)^K$.
Hence $\{\varphi_\gamma(a), \varphi_\gamma(b)\}=0$,
$\{\varphi_\beta(a), \varphi_\beta(b)\}=0$
for any
$a,b\in S(\gt g/\gt k)^K$. The homomorphisms $\varphi_\gamma$ and
$\varphi_\beta$
are surjective, so the Poisson algebras $S(\gt l_\gamma/\gt
k_\gamma)^{K_\gamma}$ and $S(\gt n)^{K_\beta}$ are commutative.
\end{proof}

\par\medskip

\begin{center}
{\large \bf 2. Properties of commutative spaces}
\end{center}

Suppose $X=G/K=(N\leftthreetimes L)/K$ is a commutative
homogeneous space. 
Denote by $P$ the ineffective kernel
of the action $L:\gt n$.
Note that $P$ is a normal subgroup of $L$ and $G$.
Due to (i) we have
$L/P\subset \Opp(\gt n)$. Hence, the stabiliser $L_\gamma$
is reductive for any $\gamma\in\gt n^*$ and the generic
stabiliser $L_*(\gt n)$ is also reductive. Condition
(i) holds \iff $L=L_*(\gt n)K$ or equivalently
$L/P$ is a product of $L_*(\gt n)/P$ and $K/(K\cap P)$. 
All nontrivial factorisations of
compact groups as products of two subgroups are
classified by Onishchik \cite{On}.
The group $L_\gamma$ is reductive, hence
the homogeneous space $L_\gamma/K_\gamma$  considered in (ii),
is commutative if and only if it is spherical.

\begin{df} Let $M$, $F$, $G$, $K$ be Lie groups, with
$F\subset M$ and $K\subset G$. The pair $(M, F)$
is called an {\it extension} of $(G,K)$ if
$$
G\subsetneqq M, \enskip M=GF, \enskip K=F\cap G.
$$
\end{df}
\noindent
Condition (i) means that $(L,K)$ is an extension of
$(L_*(\gt n), K_*(\gt n))$.

Evidently, $G/P=N\leftthreetimes (L/P)$ and
$(G/P)/[K/(K\cap P)]$ is a commutative homogeneous space of $G/P$.
In this section we are interested in commutative spaces 
satisfying 
condition 
\vskip0.5ex
\hbox to \textwidth{\qquad \quad $({\bf \ast})$\hfil $L\ne
K$ and
the action $L:\gt n$ is locally effective, i.e.,
$P$ is finite.
\hfil }
\noindent
In particular, this condition means that $L$ is compact.

\begin{lm}\label{2}
Let a symmetric pair $(M=F\times F, F)$
with a simple compact group $F$
be an extension of a spherical pair $(G, H)$.
Then $G$ contains either $F\times\{e\}$ or
$\{e\}\times F$.
\end{lm}

\begin{proof} Let $G_1$ and $G_2$ be
the images of the projections
of $G$ onto the first and the second
factors respectively. The group $G_1\times G_2$
acts spherically on $F\cong M/F\cong G/H$.
If neither $G_1$ nor
$G_2$ equals $F$, then  due to \cite[Theorem 4]{ap} we have
$\dim B(G_i(\cp))\le\dim U(F(\cp))$.
Hence, $\dim B((G_1\times G_2)(\cp))\le 2\dim U(F(\cp))<\dim F(\cp)$ and
the action $(G_1\times G_2):F$ cannot be spherical.
Assume that $G_1=F$ but $F\times\{e\}$ is not contained in
$G$. Then $G\cong F$ and $H=\{e\}$. But
the pair $(F, \{e\})$ cannot be
spherical.
\end{proof}

\begin{lm}\label{3}
Suppose a compact group $F\subset\Sp_n$ acts irreducibly
on $\qv^n$ and $F|_{\xi\qv} =\Sp_1$ for every $\xi\in\qv^n$, $\xi\ne 0$. Then
$F=\Sp_n$.
\end{lm}

\begin{proof} Let $F(\cp)\subset\Sp_{2n}(\cp)$ be
the complexification of $F$.
Then $F(\cp)$ acts on a generic subspace $\cp^2\subset\cp^{2n}$
as $\SL_2(\cp)$. Hence it acts on $\cp^{2n}$ locally transitively.
It was proved by Panyushev \cite{dp2} in a classification-free way,
that $F(\cp)=\Sp_{2n}(\cp)$.
\end{proof}

\begin{lm}\label{4}
Suppose $\gt l\subset\gt{so}(V)$ is a Lie algebra.
Let $\gt l_1$ be a  non-abelian simple ideal of $\gt l$.
Denote by $\pi$ the projection onto
$\gt l_1$. If $\pi(\gt l_*(V))=\gt l_1$ and $W_1$ is a
non-trivial  irreducible
$\gt l$-submodule of $V$ that is also non-trivial as an
$\gt l_1$-module, then
$\gt l_1=\gt{su}_2$ and $W_1$
is of the form $\mathbb H^1\otimes_{\mathbb H}\mathbb H^n$,
where $\gt l$ acts on $\qv^n$ as $\gt{sp}_n$.
\end{lm}

\begin{proof}
Set $\gt l=\gt l_1\oplus\gt l_2$.
We may assume that $V=W_1$.
The vector space $V$ can be decomposed into a tensor
product $V=V_{1,1}\otimes_{\mathbb D} V^1_1$
of $\gt l_1$ and $\gt l_2$-modules, 
where $\mathbb D$ is
one of $\rl$, $\cp$ or $\qv$. Here $\gt l_1$ acts trivially on $V_1^1$ and
$\gt l_2$ acts trivially on $V_{1,1}$. 
Both actions $\gt l_1:V_{1,1}$ and $\gt l_2:V_1^1$
are irreducible. 

Let $x=x_{1,1}\otimes x_1^1\in V$ be a non-zero
decomposable vector.
Because $V_{1,1}$ is a non-trivial irreducible 
$\gt l_1$-module, $(\gt l_1)_x\ne\gt l_1$. 
We have $\gt l_*\subset \gt l_x$ up to conjugation.
Evidently, $\gt l_x\subset\gt n_1(x)\oplus\gt n_2(x)$, where
$\gt n_i(x)=\{\xi\in \gt l_i : \xi x\in\mathbb D x\}$.
Since $\gt l_1=\pi(\gt l_*)\subset\gt n_1(x)$, 
we have $\gt n_1(x)=\gt l_1$. Hence,   
$\mathbb Dx_{1,1}$ is an $\gt l_1$-invariant 
subspace of $V_{1,1}$. Thus  $V_{1,1}=\mathbb Dx_{1,1}$
and $\gt l_1\subset\gt{gl}_1(\mathbb D)$.
If ${\mathbb D}$ equals $\rl$ or $\cp$, 
$\gt{gl}_1(\mathbb D)$ is commutative.
If ${\mathbb D}=\qv$ then $\gt l_1(x)\subset\gt{sp}_1$. 
Thus we have shown that $\gt l_1=\gt{sp}_1=\gt{su}_2$
and $W_1=\mathbb H^1\otimes_{\mathbb H}\mathbb H^n$. 
Moreover, $\gt l_2|_{\qv x}=\gt{sp}_1$.
To conclude, notice that $\gt l$ has to act
on $\qv^n$ as $\gt{sp}_n$
by Lemma \ref{3}.
\end{proof}

Set $L_*:=L_*(\gt n)$, $K_*:=K_*(\gt n)$.
Recall that there is a factorisation $L=Z(L)\times L_1\times ...\times L_m$.
Denote by $\pi_i$ the projection onto $L_i$.

\begin{thm} Suppose a commutative homogeneous space
$X=G/K$ satisfies condition  $({\bf \ast})$.
Then any non-commutative normal subgroup of $K$ distinct from
$\SU_2$ 
is contained in some simple factor of $L$.
\end{thm}

\begin{proof} Assume  $K_1$ is a normal subgroup of $K$ having
non-trivial projections onto, say, $L_1$ and $L_2$.
Consider the subgroup
$M=Z(L)\times\pi_1(K)\times\pi_2(K)\times L_3\times ...\times L_m$.
Evidently, $K\subset M$. Without loss of
generality we can replace $L$ by $M$ or better assume from the beginning that
$L_i=\pi_i(K)=\pi_i(K_1)\cong K_1$ ($i=1,2$).
Denote by $\pi_{1,2}$ the projection onto
$L_1\times L_2$. According to condition (i),
$L_1\times L_2=K_1\pi_{1,2}(L_*)$. 
Recall that due to condition (ii), $L_*/K_*$
is commutative, hence $(L_*, K_*)$ is a spherical pair. The pair
$(\pi_{1,2}(L_*), \pi_{1,2}(K_*))$ is
also spherical as a factor of a spherical pair.
Clearly,
$\pi_{1,2}(K_*)\subset\pi_{1,2}(K)\cap \pi_{1,2}(L_*)$.
Thus the symmetric pair $(L_1\times L_2, K_1)$ is an
extension of the spherical pair $(\pi_{1,2}(L_*),
\pi_{1,2}(K)\cap\pi_{1,2}(L_*))$. By Lemma \ref{2}
the group $\pi_{1,2}(L_*)$ contains
$L_1$ or $L_2$ (we can assume that
it contains $L_1$). Then $\pi_1(L_*)=L_1$
and by Lemma \ref{4} we have
$L_1=\SU_2$.
\end{proof}

In Table 1 we present the  list
of all factorisations of compact simple Lie algebras
obtained in \cite{On}.
Here $\gt g^1$, $\gt g^2$ are subalgebras of $\gt g$, $\gt g
=\gt g^1+\gt g^2$, $\gt u=\gt g^1\cap\gt g^2$.
 In all cases
$n>1$, $\varphi^1$ and $\varphi^2$ are the restrictions of
the defining representation of the complexification
$\gt g(\cp)$ to $\gt g^1(\cp)$ and $\gt g^2(\cp)$
(whose highest weights are indicated), $\varpi_m$ are the
fundamental weights, $\mymathbb $ is the trivial
representation.

\par\smallskip
\begin{center}

Table 1.
\nopagebreak
\par\smallskip

\begin{tabular}{|>{$}c<{$}|>{$}c<{$}|>{$}c<{$}|>{$}c<{$}|>{$}c<{$}|>{$}c<{$}|}
\hline
\gt g & \gt g^1 & \varphi^1 & \gt g^2 & \varphi^2 & \gt u \\
\hline
\gt{su}_{2n} & \gt{sp}_{n} & \varpi_1 &
       \gt{su}_{2n-1} & \varpi_1+\mymathbb & \gt{sp}_{n-1} \\
\cline{4-4} \cline{6-6}
 & & & \gt{su}_{2n-1}\oplus\mathbb R & & \gt{sp}_{n-1}\oplus\mathbb R \\
\hline
 \gt{so}_{2n+4} & \gt{so}_{2n+3} & \varpi_1+\mymathbb
& \gt{su}_{n+2} & \varpi_1+\varpi_{n+1} & \gt{su}_{n+1} \\
\cline{4-4} \cline{6-6}
& & & \gt{su}_{n+2}\oplus\mathbb R & & \gt{su}_{n+1}\oplus\mathbb R \\
\hline
\gt{so}_{4n} & \gt{so}_{4n-1} & \varpi_1+\mymathbb
&
\gt{sp}_{n} & \varpi_1+\varpi_1 & \gt{sp}_{n-1} \\
\cline{4-4} \cline{6-6}
& & & \gt{sp}_{n}\oplus\mathbb R & & \gt{sp}_{n-1}\oplus\mathbb R \\
\cline{4-4} \cline{6-6}
& & & \gt{sp}_{n}\oplus\gt{su}_2 & & \gt{sp}_{n-1}\oplus\gt{su}_2 \\
\hline
\gt{so}_{16} & \gt{so}_{15} & \varpi_1+\mymathbb & \gt{so}_9 & \varpi_4
                                             & \gt{so}_7 \\
\hline
\gt{so}_8 & \gt{so}_7 & \varpi_3 & \gt{so}_7 & \varpi_1+\mymathbb & G_2 \\
\hline
 \gt{so}_7 & G_2 & \varpi_1 & \gt{so}_5 &
\varpi_1+\mymathbb +\mymathbb
                     & \gt{su}_2 \\
\cline{4-4} \cline{6-6}
 & & & \gt{so}_5\oplus\mathbb R & & \gt{su}_2\oplus\mathbb R \\
\cline{4-6}
 & & & \gt{so}_6 & \varpi_1+\mymathbb & \gt{su}_3 \\
\hline
\end{tabular}
\end{center}
\par\medskip

Note that
 in Table 1 all algebras $\gt g^1$ are simple; if
$\gt g^2$ is not simple, then   $\gt g^2=\gt g^3\oplus \gt a$,
where $\gt g^3$ is simple, $\gt a=\rl$ or $\gt
a=\gt{su}_2$ and $\gt g=\gt g^1+\gt g^3$.

Suppose $X=(N\leftthreetimes L)/K$ is commutative
and $K\ne L$.  
Let $L_1\ne\SU_2$ be a simple non-commutative normal subgroup of 
$L$ such that $\pi_1(K)\ne L_1$ and $L_1\not\subset P$.
Due to Lemma \ref{4} and Theorem 1
there is a non-trivial factorisation $L_1=\pi_1(K)\pi_1(L_*)$.
Because $L_1$ is compact, 
this equality holds 
\iff
$\gt l_1=\pi_1(\gt k)+\pi_1(\gt l_*)$.
Note that if $X$ satisfies condition $({\bf \ast})$, then
no simple factor $L_i$ of $L$ is contained in $P$. 
Besides, if $X$ satisfies condition $({\bf \ast})$ and
$L_i\ne\SU_2$, then,  
by Theorem 2, $\pi_i(K)\ne L_i$ \iff $L_i\not\subset K$.  

Let $X=(N\leftthreetimes L)/K$ be
commutative and $L_1$ be a simple normal subgroup of
$L$ such that $\pi_1(K)\ne L_1$ and $L_1\not\subset P$. 
 Let $\hat V$ stand for the sum of those irreducible $L$-invariant
subspaces of $\gt n$, on which
$L_1$ acts non-trivially. Denote by $\hat P$
the identity
component of the ineffective kernel 
of $L:\hat V$. Set $\hat L=L/\hat P$.
The group $\hat L$ can be considered as the maximal connected subgroup of $L$,
which acts on
$\hat V$ locally effectively. Denote by $\hat K$ the image
of the projection
of $K$ onto $\hat L$, i.e., $\hat K\cong K/(K\cap\hat P)$.  
Clearly, $\hat V$ need not be
a subalgebra of $\gt n$, but it
can be considered as a factor algebra. Let $\hat{\gt a}$ stand for the maximal
$L$-invariant subspace of $\gt n$ on which $L_1$ acts trivially.
Evidently, $\hat{\gt a}$ is a subalgebra. Moreover, because different
$L$-invariant summands
of $\gt n$ commute, see \cite[Proposition 15]{Vin}, it is an ideal.
We identify $\hat V$ with $\gt n/\hat{\gt a}$.

\begin{prop} All triples $(\hat L, \hat K, \hat V)$ which can be
obtained from a commutative space $X$ as a result of the
procedure described above
are contained in Table 2b.
\end{prop}

\begin{proof} The homogeneous space $(\hat V\leftthreetimes L)/K$
is commutative.
Let $L=L_1\times L^1\times Z(L)$,
where $L^1$ is the product of all simple factors
of $L$ except $L_1$. The space $\hat V$ can be
represented as a sum
$$
\hat V=(V_{1,1}\otimes_{\mathbb D_1} V^1_1)\oplus\ldots \oplus
     (V_{1,p}\otimes_{\mathbb D_p} V^1_p),
$$
where $V_{1,i}$ are irreducible $L_1$-modules, $V_i^1$ are $L^1$-modules,
$L_1$ acts trivially on each $V_i^1$ and $L^1$ acts trivially
on each $V_{1,i}$.
In each summand the tensor product is taken over the skew-field
$\mathbb D_i$, which equals
$\rl$, $\cp$ or $\qv$ depending on
$V_{1,i}$ and $V_i^1$.

Set $V_1:=V_{1,1}\oplus \ldots \oplus
 V_{1,p}$. Let  $x_1=x_{1,1}+\ldots +x_{1,p}$ be
a generic vector of $V_1$. Then there is a sum of decomposable vectors
$x:=x_{1,1}\otimes x^1_1+\ldots +x_{1,p}\otimes x^1_p\in V$,
such that $x^1_i$ are linearly independent.

Note that $(\pi_1(L_x), \pi_1(K_x))$ is a spherical pair as
a factor of the spherical pair $(L_x, K_x)$.
Clearly,
$\pi_1(K_x)\subset\pi_1(L_x)\cap\pi_1(K)$.
According to condition (i), $L_1=\pi_1(L_x)\pi_1(K)$.
Hence, the pair $(L_1, \pi_1(K))$ is an extension of
the spherical pair $(\pi_1(L_x), \pi_1(L_x)\cap\pi_1(K))$.

The generic stabiliser
$(L_1)_*(V_1)$ is defined up to  conjugation.
We may assume that $(L_1)_*(V_1)=(L_1)_{x_1}\subset\pi_1(L_x)$.
Obviously,
$\pi_1(L_x)x_1\subset\mathbb D_1 x_{1,1}+\ldots+\mathbb D_p x_{1,p}$.
Hence, $(L_1)_*(V_1)$ is a normal subgroup of
$\pi_1(L_x)$ and
$\pi_1(L_x)/(L_1)_*(V_1)$ is locally isomorphic to a
direct product of
$(\SU_2)^k$ and $(\Upp_1)^n$.
Recall that there is a non-trivial factorisation
$L_1=\pi_1(L_x)\pi_1(K)$. Looking in 
Table~1 one can see that
$(L_1)_*(V_1)^0\ne\{E\}$, $L_1=\pi_1(L_x)(L_1)_*(V_1)$,
and $\pi_1(L_x)/(L_1)_*(V_1)$ is locally isomorphic to
$\SU_2$ or $\Upp_1$, or  trivial. To conclude the proof,
we need the following lemma.

\begin{lm}\label{5} All triples $(L_1, \pi_1(K), V_1)$ which
can be obtained as a result of
the above procedure
are contained in Table 2a.
\end{lm}

\begin{proof} As we already know,
$(L_1, \pi_1(K))$ is an extension of a spherical pair
 $(\pi_1(L_x),
 \pi_1(L_x)\cap\pi_1(K))$; in particular, it is spherical.
The Lie algebra $\gt l_1$ is simple, hence this extension is 
contained in Table 1. 
Suppose $\gt l_1=\gt g$ for some $\gt g$ from Table 1
and $(\gt g, \gt g_i)$ is an extension
of a spherical pair $(\gt g_j, \gt u)$. Then $(\gt
l_1)_*(V_1)$ equals to one of the following three algebras:
$\gt g_j$, $\gt g_j'$, $\gt{g}_{j}/\gt{sp}_1$. The last case
is only possible if $\gt g_j=\gt{sp}_{n}\oplus\gt{sp}_1$. We
have to check if any of these three algebras is a generic
stabiliser for some representation of $\gt l_1$. The
representations of complex simple algebras with non-trivial
generic stabiliser are classified by \'Elashvili \cite{al}.
We are interested in the real forms of 
the orthogonal representations with non-trivial
generic stabiliser. 

The algebra $\gt l_1$ is either $\gt{su}_{2n}$ or $\gt{so}_m$. 
If $\gt l_1=\gt{su}_{2n}$, $(\gt
l_1)_*(V_1)$ have to be one of the following three algebras 
$\gt{sp}_n$, $\gt{su}_{2n-1}$, $\gt{u}_{2n-1}$. According to 
\cite{al}, $V_1$ is either $\cp^{2n}$ or $\rl^6$ for 
$\gt l_1=\gt{su}_4$. 

Suppose $\gt l_1=\gt{so}_m$. According to \cite{al}, 
if $m>15$, then $V_1$ is the sum of $k$ copies 
of $\rl^m$ and $(\gt l_1)_*(V_1)=\gt{so}_{m-k}$. 
Table 1 and Kr\"amer's classification \cite{kr}
tell us that $k=1$, so $\gt g_j=\gt{so}_{m-1}$ and 
$(\gt g_j,\gt u)$ is spherical only in one case, 
namely $(\gt{so}_{2n+3}, \gt{su}_{n+1}\oplus\rl)$. 
For smaller $m$ one has to check several cases by direct 
computations. The result is given in rows 2a, 2b, 4a
and 4b of Table 2a.   

We assume that 
a triple $(L_1, \pi_1(K), W)$ is contained in Table 2a, if 
$W$ is an $L_1$-invariant subspace of some $V_1$ and  
$(L_1, \pi_1(K), V_1)$ is a triple precisely listed in Table~2a. 
\end{proof}

Proceed to the proof of Proposition 2. As can be easily seen,
Table 2a does not contain
symplectic representations, i.e., those that could be a factor
in a tensor product over $\qv$.
It contains only one Hermitian representation (in the first line).

Set $W_1:=V_{1,1}\otimes V^1_1$. Our next goal is to describe all
possible $W_1$. The pair $(L_1, \pi_1(K), V_{1,1})$ is contained
in Table 2a, i.e., $V_{1,1}$ is an $L_1$-invariant
subspace of $V_1$. In particular, the tensor product in $W_1$
is taken over $\rl$ or $\cp$.

Set $n_1=\dim V_{1,1}$, $n_2=\dim V_1^1$ and $s=min(n_1,
 n_2)$ (here we consider the dimensions over the same
field as the tensor product in $W_1$.)
Let us prove that up to conjugation
 $$
 [\pi_1(L_*(W_1))]'=[(L_1)_{\xi_1}\cap ...\cap
 (L_1)_{\xi_s}]',
 $$
where $\xi_i\in V_{1,1}$ are linear independent
generic vectors. In particular, this equality means that
if $\dim V_{1,1}\le\dim V^1_1$, then $[\pi_1(L_*(W_1))]'=\{E\}$.
Note that the right hand side is  evidently a subset of the
left hand side. So, it is enough to
proof the inclusion ``$\subset$''.

The group $L_1$ is compact, so we can assume that $L_1\subset\Opp_{n_1}$
if $\mathbb D_1=\rl$ or  $L_1\subset \Upp_{n_1}$ if
$\mathbb D_1=\cp$.
Without loss of generality it can be also assumed that
the considered action is $L_1\cdot\Opp_{n_2}:V_{1,1}\otimes_{\rl} V^1_1$ or
 $L_1\cdot\Upp_{n_2}:V_{1,1}\otimes_{\cp} V^1_1$. The left hand side
of the proving equality
could become only larger after such replacement. The required
inclusion can be deduced from
a well known fact: if $n\le m$, then  the generic stabilisers for
$\Opp_n\times\Opp_m:\rl^n\otimes\rl^m$
and $\Upp_n\times\Upp_m:\cp^n\otimes\cp^m$ are
 $\Opp_{m-n}\times(\mathbb Z/2\mathbb Z)^n$ and $\Upp_{m-n}\times
 \Upp_1^n$, respectively.

Recall that there should be factorisations $L_1=\pi_1(\hat
 K)\pi_1(\hat L_*)$ and  $L_1=\pi_1(\hat K)[\pi_1(\hat L_*)]'$.
This imposes rather strong restrictions on $s$. Thus
$s=1$ in cases 1a, 1b, 2b, 3 and 4b; $s<3$ in case 2a; $s<4$ in case 4a.
If $V_1$ is irreducible $L_1$-module, then $\hat V=W_1$. 
Hence, $\hat V=V_1$ in cases 1a, 1b, 2b, 3 and 4b.  
Note that in general
$$ 
 \pi_1(\hat L_*(\hat V))\subset \pi_1(\hat L_*(V_{1,1}\otimes
 V^1_1)\cap\ldots\cap\hat L_*(V_{1,p}\otimes V^1_p)).
$$
So there are only a few possibilities for $\hat V$ and hence for 
$(\hat L, \hat K)$.
Namely, $p\le 2$, $(\dim V_1^1+\dim V_2^1)\le 2$ in case 2a; 
and $p\le 3$, $(\dim V_1^1+\dim V_2^1+\dim V_3^1)\le 3$ in case 4a. 
For each representation one can easily verify whether conditions 
(i) and (ii) hold.
The result is shown in Table 2b.
\end{proof}

\begin{rmk} In the first row of Table 2b $\hat V$ is written
as $\cp^{2n}(\oplus\rl)$. Of course, $\rl$ is a trivial $L$-module.
So, formally it cannot be put into the table. But it is done here to indicate
that in this case $\hat V$ can be a non-commutative subspace of $\gt n$.
\end{rmk}

\begin{center}
 \par\medskip
 \hbox to 17cm {
 \begin{tabular}{|c|>{$}c<{$}|>{$}c<{$}|>{$}c<{$}|}
 \multicolumn{4}{c}{Table 2a.} \\
 \hline
  & L_1 & \pi_1(K) & V_1 \\
 \hline
 1a & \SU_{2n} & \Sp_n & \cp^{2n} \\
 \hline
 1b & \SU_4 & \Upp_3 & \rl^6 \\
 \hline
 2a & \SO_7 & G_2 & \rl^7\oplus\rl^7 \\
 \hline
 2b & \Spin_7 & \Spin_6 & \rl^8 \\
    &   & \Spin_5\cdot\Upp_1 & \\
 \hline
 3 & \SO_{2n} & \SU_n & \rl^{2n} \\
   &  & \Upp_n & \\
 \hline
 4a & \SO_8 & \Spin_7 & \rl^8\oplus\rl^8\oplus\rl^8 \\
 \hline
 4b & \SO_8 & \Sp_2\times\SU_2 & \rl^8 \\
 \hline
 \end{tabular} \hfill
 \begin{tabular}{|c|>{$}c<{$}|>{$}c<{$}|>{$}c<{$}|}
 \multicolumn{4}{c}{Table 2b.} \\
 \hline
  & \hat L & \hat K & \hat V \\
 \hline
 1a & {\rm (S)U}_{2n} & \Sp_n(\cdot\Upp_1) & \cp^{2n}(\oplus\rl) \\
 \hline
 1b & \SU_4 & \Upp_3 & \rl^6 \\
 \hline
 2a & \SO_7 & G_2 & \rl^7 \\
 \hline
 2b & \Spin_7 & \Spin_6 & \rl^8 \\
 \hline
 3 & \SO_{2n} & \Upp_n & \rl^{2n} \\
 \hline
 4a & \SO_8\times\SO_2 & \Spin_7\times\SO_2 & \rl^8\otimes\rl^2 \\
 \hline
 4b & \SO_8 & \Spin_7 & \rl^8\otimes\rl^2 \\
 \hline
 4c & \SO_8 & \Spin_7 & \rl^8 \\
 \hline
 4d & \SO_8 & \Sp_2\times\SU_2 & \rl^8 \\
 \hline
 \end{tabular} }
 \end{center}
 \par\medskip

So far nothing was said about the Lie algebra structure on $\hat V$. Denote by
$\hat{\gt n}$ the subalgebra of $\gt n$ generated by $\hat V$, i.e.,
$\hat{\gt n}=\hat V+[\hat V, \hat V]$.

Suppose $[\hat V, \hat V]=0$.  Then we have a commutative
homogeneous space $(\hat N\leftthreetimes\hat L)/\hat K$,
where $\hat N$ is a simply connected commutative Lie group.

Suppose $X=(N\leftthreetimes L)/K$ 
is commutative, $L_1\ne\SU_2$ is a simple normal subgroup of 
$L$, such that $\pi_1(K)\ne L_1$, $L_1\not\subset P$ 
and a triple $(\hat L, \hat K, \hat V)$ corresponds to 
$L_1$ in the aforementioned way. By Proposition 2, 
$(\hat L, \hat K, \hat V)$ is contained in Table 2b.  

\begin{lm}\label{6} If $(\hat L, \hat K, \hat V)$ 
is a triple from Table 2b but
not from the first row of it, 
then $\hat V$ has to be a commutative subspace of $\gt n$.
\end{lm}

\begin{proof}
Assume that $[\hat V, \hat V]\ne 0$.
There in an $L$-invariant surjection 
$\Lambda^2\hat V\mapsto[\hat V,\hat V]$. 
Because representation of $L$ in $\Lambda^2\hat V$ 
is completely reducible, 
$[\hat V, \hat V]$ can be regarded as an 
$L$-invariant subspace of $\Lambda^2\hat V$. 
If $[\hat V, \hat V]$ is a non-trivial $L_1$-module, then 
it is an $\hat L$-invariant subspace of 
$\hat V$ (by definition of $\hat V$). 
Recall that $\gt n$ is nilpotent.
Hence, $[\hat V, \hat V]\ne\hat V$. 
In particular, if $\hat V$ is 
an irreducible $\hat L$-module, 
then $L_1$ has to act on $[\hat V, \hat V]$ trivially.
In all rows of Table 2b except 4b the representation
$\hat L:\hat V$ is irreducible. 
The space $\Lambda^2\hat V$ contains non-trivial $\hat L$-invariants only
in cases  1a, 4a and 4b.

Consider case 4a. Here $\hat K_*(\gt m)=\SU_4\times\SO_2$ and
$\rl^8{\otimes}\rl^2=\cp^4{\otimes_\rl}\cp\cong\cp^4{\oplus}\cp^4$
is a sum of two isomorphic $K_*(\gt m)$-modules.
According to \cite[Proposition 15]{Vin}, these two submodules 
are commutative and  
commute with each other. Hence, $[\hat V, \hat V]=0$. 
In case 4b $\hat V=\rl^8\oplus\rl^8$ is a sum 
of two isomorpich $\hat L$-modules. Hence, it is 
commutative. 
\end{proof}

The first case is different. It corresponds
to six commutative spaces, namely $\hat L$ can be 
either $\SU_{2n}$ or $\Upp_{2n}$.
In the second case there are also two possibilities $\hat K=\Sp_n$ or
$\hat K=\Sp_n\times\Upp_1$.
Independently $\hat{\gt n}$ can be either
$\cp^{2n}$ or $\cp^{2n}\oplus\rl$, with
$\hat N$ being commutative or
the Heisenberg group $H_{2n}$.
So, the first row gives rise to six commutative spaces and
each of the other 
to  only one. In each case the commutativity can easily be
proved by means of Theorem 1.
Two cases are considered in detail.

\begin{ex} Let us prove that
 $(H_n\leftthreetimes \Upp_{2n})/\Sp_n$ is
commutative. Since $\Upp_{2n}$ and $\Sp_n$
are transitive on the sphere in $\cp^{2n}$,
$\rl[\cp^{2n}]^{\Upp_{2n}}=\rl[q]=\rl[\cp^{2n}]^
{\Sp_n}$, where $q$ is the  invariant of degree 2.

The generic stabiliser for $\Sp_n:\cp^{2n}$ is equal to
$\Sp_{n-1}$. The space
$\Upp_{2n-1}/\Sp_{n-1}$ is a compact real form of the complex
spherical space $\GL_{2n-1}(\cp)/\Sp_{2n-2}(\cp)$, and, hence, is
commutative.

So, only (iii) is left. Here we have $\gt
m=\gt{u}_{2n}/\gt{sp}_n=\bigwedge^2\cp^{2n}$. It is a classical
result that
$K_*(\bigwedge\cp^{2n})=\underbrace{\SU_2\times...\times\SU_2}_{n}$.
As a $K_*(\gt m)$-module $\gt
n=\gt v_1\oplus...\oplus\gt v_n\oplus\rl$, where $\gt v_i=\cp^2$
for every $i$. Each $\gt v_i$ is acted upon by its own $\SU_2$.
Note that $[\gt v_i, \gt v_j]=0$ for $i\ne j$. For $K_*(\gt
m)$-invariants we have $S(\gt
n)^{K_*(\gt m)}=\rl[t_1,...t_n,\xi]$, where $t_i$ is the quadratic
$\SU_2$ invariant in $S^2(\gt v_i)$ and $\xi\in\gt h_n'$.
Evidently, $t_i$ and  $t_j$
commute as elements of the  Poisson algebra $S(\gt n)$, and $\xi$
lies in the centre of $S(\gt n)$.
\end{ex}

\begin{ex}
The homogeneous space $(\rl_{2n}\leftthreetimes \SO_{2n})/\Upp_n$ is
also commutative. Here $\gt n$ is commutative, so we do not need to check
condition (iii). For (i) we have
$\rl[\rl^{2n}]^{\SO_{2n}}=\rl[q]=\rl[\rl^{2n}]^{\Upp_n}$.
It can be easily seen that $L_*=\SO_{2n-1}$ and $K_*=\Upp_{n-1}$.
The corresponding homogeneous space  $\SO_{2n-1}/\Upp_{n-1}$ is spherical
by the Kr\"amer's classification \cite{kr}.
\end{ex}

Note that the triple $(\SU_4, \U_3, \rl^6)$ is locally
isomorphic to a triple from 
row 3 of Table 2b with n=3. 

\begin{prop}\label{simple} Suppose $X=(N\leftthreetimes L)/K$ is 
an indecomposable commutative space, 
$\gt n\ne 0$, 
$L$ is simple and $L\ne K$. Then $X$ is 
a homogeneous space from Table~2b.   
\end{prop}

\begin{proof} 
The action $L:\gt n$ is non-trivial, otherwise 
$X=N\times(L/K)$. Consider the 
triple $(L, K, \hat V)$ corresponding to $L_1=L$. 
By Proposition 2, it is contained in Table 2b. 
Assume that $\hat{\gt n}\ne\gt n$ and 
let $\gt a$ be an $L$-invariant complement of 
$\hat{\gt n}$ in $\gt n$. By definition of $\hat V$, 
$L$ acts on $\gt a$ trivially. Thus $\gt a$ is an abelian ideal
and $X$ is decomposable $X=((\hat N\leftthreetimes L)/K)\times A$, 
where $A\subset N$ and $\Lie A=\gt a$.  
\end{proof}

\begin{lm}\label{7} Let $(N\leftthreetimes L)/K$ be a commutative
homogeneous space satisfying
condition $({\bf \ast})$ and
$K_1\cong\SU_2$ a normal subgroup of $K$. Then either
$K_1\subset\SO_8$, where $\SO_8$
is a simple direct factor of $L$ contained in the row
4d of Table 2b, or $K_1$ is the diagonal of
a product of at most three simple direct 
factors of $L$ isomorphic to $\SU_2$.
\end{lm}

\begin{proof} Suppose $\pi_i(K_1)\ne\{E\}$
and $L_i\ne\SU_2$. Then $\pi_i(K)\ne L_i$. 
Consider a triple $(\hat L, \hat K, \hat V)$ 
corresponding to $L_i$. By Proposition 2, it is 
contained in Table 2b. Note that $K_1$ is a normal subgroup 
of $\hat K$. Thus $\hat K=\Sp_2\times\SU_2$,
$\hat L=\SO_8$. Assume that $K_1$ has a non-trivial projections onto
some other simple direct factor of $L$. 
Then the pair $(\SO_8\times\SU_2, \Sp_2\times\SU_2)$
with $\SU_2$ embedded into $\SO_8$ as a centraliser of
$\Sp_2$ and into $\SU_2$ isomorphically,
should be spherical, but it is not.
To conclude, note that the pair
$(\SU_2\times\SU_2\times\SU_2\times\SU_2, \SU_2)$ is not spherical either.
\end{proof}

\par\medskip


\begin{center}
{\large \bf 3. Principal commutative spaces}
\end{center}

Let $G/K=(N\leftthreetimes L)/K$ be a commutative homogeneous space. 
Recall that by our assumptions $L=Z(L)\times L_1\times\dots L_m$, 
where $Z(L)$ is the connected centre of $L$. The 
ineffective kernel of $L:\gt n$ is denoted by $P$.
Denote by $Z(K)$ the connected centre of $K$.
       Decompose $\gt n/\gt n'$ into a sum of irreducible
$L$-invariant subspaces
$\gt n/\gt n'=\gt w_1\oplus ... \oplus\gt w_p$.

\begin{df} Let us call a commutative homogeneous space $G/K$
{\it principal} if $P$ is semisimple,
 $Z(K)=Z=Z(L)\times (L_1\cap Z)\times...\times(L_m\cap Z)$
 and
$Z(L)=C_1\times...\times C_p$, where 
$C_i\subset\GL(\gt w_i)$.
\end{df}

The classification of commutative homogeneous spaces can
be divided in two parts: the classification of
principal commutative spaces and description of
the possible centres of $L$ and $K$ in the general case.
Note that Table 2b contained only two non-principal 
homogeneous spaces, namely 
$(H_{2n}\leftthreetimes\U_{2n})/\Sp_n$ and  
$(\cp^{2n}\leftthreetimes\U_{2n})/\Sp_n$.

\begin{ex}
We have proved that the homogeneous spaces
$X=(H_{2n}\leftthreetimes\SU_{2n})/\Sp_n$ and
$Y=(H_{2n}\leftthreetimes\Upp_{2n})/\Sp_n$, where
$N=H_{2n}$ is the Heisenberg group with
$\Lie H_{2n}=\gt h_{2n}=\cp^{2n}\oplus\rl$, are commutative.
Denote by $X^m=(H_{2n}\leftthreetimes\SU_{2n})^m/(\Sp_n)^m$
the product of $m$ copies of $X$. Let $L$ and $K$ be the products of
$m$ copies of $\SU_{2n}$ and $\Sp_n$ respectively.
Suppose $F_1\subset F_2\subset(\Upp_1)^m$. Then, evidently,
$((H_{2n})^m\leftthreetimes(F_2\times L))/(F_1\times K)$
is a commutative homogeneous
space, in general indecomposable but non-principal.
\end{ex}

The above example is rather simple, because there are no conditions
on $F_1$ and $F_2$.
In other cases the situation is more difficult.
For commutative homogeneous spaces of reductive Lie groups
the description of possible centres of $L$ and $K$ is given
in \cite{y2}. The same problem for
commutative homogeneous spaces
$(H_n\leftthreetimes K)/K$, where $K\subset\Upp_n$
is solved in \cite{L} and \cite{be-ra2}. 
In the present article  we concentrate on principal commutative spaces.

\begin{thm} Let $X=(N\leftthreetimes L)/K$ be a maximal 
indecomposable
principal commutative
homogeneous space satisfying
condition $({\bf \ast})$. Then either 
$X$ is contained in Table 2b (and $L'$ is simple); 
or $(L, K)$ is
isomorphic to a product of pairs 
$(\SU_2\times\SU_2\times\SU_2, \SU_2)$,
$(\SU_2\times\SU_2, \SU_2)$ or $(\SU_2, \Upp_1)$ and a pair
$(K^1, K^1)$, where $K^1$ is a compact Lie group. 
\end{thm}

\begin{proof} Suppose there is a simple normal subgroup $L_i\ne\SU_2$ 
of $L$, which is not contained in $K$. Then by Theorem 2 
$\pi_i(K)\ne L_i$. Consider the corresponding triple 
$(\hat L, \hat K, \hat{\gt n})$. It is a homogeneous space from 
Table 2b. Recall that by definition $\hat L$ is the maximal connected 
subgroup of $L$ acting on $\hat V$ locally effectively,
$\hat P$ is the identity component of the ineffective kernel 
of $L:\hat V$ and $L=\hat L\cdot\hat P$. 
Because $X$ is principal, $Z(L)=\hat C\times C^1$, where 
$\hat C=\GL(\hat V)\cap Z(L)$, hence $\hat C\subset\hat L$,
and $C^1\subset\hat P$. In particular, 
$L=\hat L\times\hat P$. Similarly, the connected centre $Z=Z(K)$ 
is a product $Z=Z(L)\times \hat Z\times Z^1$, where
$\hat Z\subset\hat L'$ and $Z^1\subset\hat P'$. 
According to Theorem 2 each normal 
subgroup $K_i\not\cong\SU_2$ of $K$ is contained 
in some simple direct factor of $L$, hence either in $\hat L$ 
or in $\hat P$. Suppose a normal subgroup
$K_j\cong\SU_2$ of $K$ is not contained in any
simple direct factor of $L$. Then, 
by Lemma~\ref{7}, it is diagonal in a product of 
at most three direct factors of $L$ isomorphic to $\SU_2$. 
The group $\hat L$ has no normal subgroups isomorphic to $\SU_2$.
Hence, $K_j\subset\hat P$.
Thus $\hat K=\hat L\cap K$. 
Moreover, $K=\hat K\times F$, where $F\subset\hat P$. 
Recall that $\gt n=\hat V\oplus\hat{\gt a}$, 
where $\hat{\gt a}$ is an ideal and $L_i$ acts on 
$\hat{\gt a}$ trivially. Let $\hat A\subset N$ be a corresponding 
connected subgroup. Note that either 
$\hat L=L_i$ or $\hat L=\U_1\times L_i$. Anyway
$\hat L$ acts on $\hat{\gt a}$ trivially. 
Thus in case $[\hat V, \hat V]=0$, i.e.,
$\hat{\gt n}=\hat V$, we have obtained a decomposition 
$X=((\hat N\leftthreetimes\hat L)/\hat K)\times
    ((\hat A\leftthreetimes\hat P)/F)$. 
But $X$ in indecomposable, hence 
$X=(\hat N\leftthreetimes\hat L)/\hat K$ and it is 
contained in table 2b. 

There is only one possibility for 
non-commutative $\hat{\gt n}$, namely 
$L_i=\SU_{2n}$, $\gt n=\hat V\oplus\gt z$, 
where $\gt z\in\hat{\gt a}$ and $\gt z\cong\rl$ is 
a trivial $L$-module. Let $\gt a$ be 
an $L$-invariant complement 
of $\gt z$ in $\hat{\gt a}$. If $\gt a$ is a
subalgebra (then it is an ideal), 
we again have a decomposition of $X$. 
Assume that $\gt z\subset[\gt a,\gt a]$. 
We have an $L$-invariant decomposition
$\gt n=\gt a\oplus\gt z\oplus\hat V$, 
where $[\gt a,\gt a]\subset\gt a\oplus\gt z$, 
$[\hat V, \hat V]=\gt z$ and $[\gt z,\gt n]=0$. 
Thus $X$ is a central reduction of 
$((\hat N\leftthreetimes\hat L)/\hat K)\times
    ((\hat A\leftthreetimes\hat P)/F)$  
by a one dimensional subgroup embedded diagonally 
into $\hat N'\times\hat A'$. Hence $X$ is not maximal. 

We have proved that if there is a simple normal subgroup $L_i\ne\SU_2$ 
of $L$, which is not contained in $K$, then 
$X$ is a homogeneous space from Table 2. 
If this is not the case, then 
the spherical pair $(L, K)$ is a product of
the ``$\SU_2$-pairs'' and
$(K^1, K^1)$, where
$K^1$ contains the connected center of $L$ and 
all its simple normal subgroups different from $\SU_2$. 
\end{proof}

The homogeneous space $(\hat N\leftthreetimes\hat L)/\hat K$
corresponding to the first row of
Table 2b is maximal if and only if $\gt n=\gt h_{2n}$, 
and it is principal \iff 
$\hat L=\SU_{2n}$, $\hat K=\Sp_n$ or 
$\hat L=\U_{2n}$, $\hat K=\U_1\cdot\Sp_n$. 
Homogeneous spaces corresponding to other rows 
of Table 2b are maximal and principal.

Let $X$ be a commutative homogeneous space. Denote by
$(L^{\triangle}, K^\triangle)$ a spherical subpair
of $(L,K)$ of the type
 $(\SU_2\times\SU_2\times\SU_2, \SU_2)$, $(\SU_2, \Upp_1)$ or
$(\SU_2\times\SU_2, \SU_2)$ and by $\pi^\triangle$
a projection onto $L^{\triangle}$.

\begin{lm}\label{8} If
$(L^{\triangle},
 K^\triangle)=(\SU_2\times\SU_2\times\SU_2, \SU_2)$ or
 $(L^{\triangle}, K^\triangle)=(\SU_2, \Upp_1)$ then
 $\pi^{\triangle}(L_*)=L^\triangle$, if $(L^{\triangle},
 K^\triangle)=(\SU_2\times\SU_2, \SU_2)$ then
 $\pi^{\triangle}(L_*)$ equals $L^\triangle$,
 or $\SU_2\times\Upp_1$.
\end{lm}

\begin{proof} The group $\SU_2$ has only trivial factorisations, besides,
$(\pi^{\triangle}(L_*), \pi^\triangle(L_*)\cap
 K^{\triangle})$ is a spherical pair. In particular,
 $\pi^\triangle(L_*)\cap K^{\triangle}$ is not empty.
This reasoning explains the second  and
the third cases. It remains to observe that in
the first case the group
 $\pi^{\triangle}(L_*)$ can not be
 $\SU_2\times\SU_2\times\Upp_1$, because the pair
$(\SU_2\times\SU_2\times\Upp_1, \Upp_1)$ is not spherical.
\end{proof}

We complete our classification modulo a description of
possible actions of normal subgroups of $L$ isomorphic
to $\SU_2$ on $\gt n$. This description
is a very intricate
problem, which may be a subject of another article.

\par\medskip

\begin{center}
 {\large \bf 4. The ineffective kernel}
 \end{center}

Suppose $X=(N\leftthreetimes L)/K$ is a commutative 
homogeneous Riemannian space.
Let $P$ be the ineffective kernel of the action $L:\gt n$. Then
$L$ can be decomposed as
$L=P\cdot L^\diamond$, where $L^\diamond$ is the maximal
connected normal subgroup of $L$
acting on $\gt n$ locally effectively.
We assume that $G$ is not reductive,
hence $P\ne L$.
From the classification of spherical pairs we know that
any normal subgroup $K_1$ of $K$ not locally isomorphic to $\SU_2$
can have non-trivial projections only on two
different simple factors of $L$.

\begin{lm}\label{9} Let $X$ be commutative.
Suppose a normal subgroup $K_1\ne\SU_2$ of $K$ is not
contained in either $P$ or $L^\diamond$. Then there are
simple factors $P_1$, $L_1^\diamond$ of $P$, $L^\diamond$
such that $K_1\subset P_1\times L_1^\diamond$, $P_1\cong L
_1^\diamond\cong K_1$. Moreover, either $K_1=\SO_{n+1}$,
where $n\ge 4$; or $K_1=\SU_{n+1}$, where $n\ge 2$.
\end{lm}

\begin{proof} It can be seen from
the classification of spherical pairs that $K_1\subset
L_i\times L_j$. We can assume that $K_1\subset P_1\times
L_1^\diamond$. The action $K_1:\gt n$ is non-trivial,
otherwise $K_1$ would be a subgroup of $P$. Denote by
$\pi_1^K$ the projection onto $K_1$ in $K$ and by
$\pi_{1,1}$ the projection onto $P_1\times L_1^\diamond$ in
$L$. By Lemma \ref{4}, $\pi_1^K(K_*)\ne K_1$. Recall that
$(L_*, K_*)$ is spherical. Hence, the pair $(\pi_{1,1}(L_*),
\pi_{1,1}(K_*))$ is also spherical. Note that $L_*=P\cdot
L^\diamond_*(\gt n)$. Hence,
$\pi_{1,1}(L_*)=P_1\times\pi_1^\diamond(L_*)$, 
where $\pi_1^\diamond$ is a projection onto $L_1^\diamond$ in
$L$.

We claim that $(K_1\times\pi_1^K(K_*), \pi_1^K(K_*))$ is
spherical. Without loss of generality, we can assume that
$P_1\cong L_1^\diamond\cong K_1$. If it is not the case, we
replace $L$ by a smaller subgroup containing $K$, namely
each of $P_1$ and $L_1^\diamond$ is replaced by a projection
of $K$ onto it. We illustrate the embedding
$\pi_{1,1}(K_*)\subset\pi_{1,1}(L_*)$ by the following
diagram.
\begin{center}
\begin{tabular}{c}
\xymatrix@R+1mm@C-6mm{
 {\pi_{1,1}(L_*)} & {\cong} &
{K_1} & \times & {\pi_1^\diamond(L_*)} \\
 {\pi_{1,1}(K_*)} & {=} & & \pi_1^K(K_*)
\ar@{-}[ul] \ar@{-}[ur] & \\ }
\end{tabular}
\end{center}
Because the pair $(\pi_{1,1}(L_*), \pi_{1,1}(K_*))$ is
spherical,
 $(K_1\times\pi_1^K(K_*), \pi_1^K(K_*))$ is also spherical.
According to the classification of spherical pairs, there
are only two possibilities: either $K_1=\SO_{n+1}$,
$\pi_1^K(K_*)=\SO_{n}$; or $K_1=\SU_{n+1}$,
$\pi_1^K(K_*)=\Upp_{n}$. Assume that $P_1\cong K_1$. 
If $L_1^\diamond$ is
larger than $K_1$, then $(P_1\times L^\diamond_1,
\pi_{1,1}(K))$ is one of the following three pairs:
$(\SO_{n+1}\times\SO_{n+2}, \SO_{n+1})$;
 $(\Sp_2\times\Sp_{m+2}, \Sp_2\times\Sp_m)$;
 $(\SU_{n+1}\times\SU_{n+2}, \U_{n+1})$.
Recall that
$L_1^\diamond=\pi_1^\diamond(L_*)\pi_1^\diamond(K)$. The
group $\Sp_{m+2}$ has no non-trivial factorisation, hence
the second case is not possible. Also, we know that
$\pi_1^\diamond(L_*)\cap\pi_1^\diamond(K)$ contains $\SO_n$
or $\U_n$, depending on $K_1$. Thus, only one possibility is
left $(K_1\times L^\diamond_1,
\pi_{1,1}(K))=(\SU_3\times\SU_4, \U_3)$. According to Table
1, $\pi_1^\diamond(L_*)=\Sp_2$,
$\Sp_2\cap\U_3=\Sp_1\times\U_1$. We have
 $\pi_{1,1}(L_*)=\SU_3\times\Sp_2$. The subgroup
$\pi_{1,1}(K_*)$ is contained in $\Sp_1\times\U_1$, which is
not spherical in $\SU_3\times\Sp_2$. Hence, in this case the
pair $(\pi_{1,1}(L_*), \pi_{1,1}(K_*))$ is not spherical. We
illustrate this case by the following diagram.
\begin{center}
\begin{tabular}{c}
\xymatrix@R+1mm@C-6mm{
 {\pi_{1,1}(L_*)} & {=} &
{\SU_3} & \times & {\Sp_2} \\
{\pi_{1,1}(K_*)} & {\subset} & & {\Sp_1\times\U_1}
\ar@{-}[ul] \ar@{-}[ur] & \\ }
\end{tabular}
\end{center}
\noindent
Thus, we have proved that $L_1^\diamond\cong K_1$.

To conclude, we show that $P_1\cong K_1$. Denote by
$\pi_1^P$ the projection onto $P_1$. We can decompose
$\pi_1^P(K)$ into a locally direct product $\pi_1^P(K)\cong
F\cdot K_1$. The subgroups $F\cdot K_1$ and
$F\cdot\pi_1^K(K_*)$ are spherical in $P_1$. Moreover, the
pairs $(P_1\times K_1, F\times K_1)$ and
$(P_1\times\pi_1^K(K_*), F\times \pi_1^K(K_*))$ are also
spherical. There are the same three possibilities for
$(P_1\times K_1, F\times K_1)$, in which $P_1\ne K_1$,
namely: $(\SU_{n+2}\times\SU_{n+1}, \Upp_{n+1})$,
$(\Sp_{m+2}\times\Sp_2, \Sp_{m}\times\Sp_2)$ and
$(\SO_{n+2}\times\SO_{n+1}, \SO_{n+1})$. But even the pair
$(P_1, F\cdot\pi_1^K(K_*))$ is not spherical in any of these
cases.
\end{proof}

\begin{ex} We show that the homogeneous spaces 
$((\rl^n\leftthreetimes\SO_n)\times\SO_n)/\SO_n$ and 
$((H_n\leftthreetimes\U_n)\times\SU_n)/\U_n$ are commutative. 
We have $L_*=\SO_n\times\SO_{n-1}$ for the first space 
and $L_*=\SU_n\times\U_{n-1}$ for the second one, 
$K_*$ is either $\SO_{n-1}$ or $\U_{n-1}$. 
The stabiliser $L_*$ contains the first direct factor
and $K$ is the diagonal multiplied by the connected centre
of $L$. Hence, $L=L_*K$. According to \cite{br} and 
\cite{m}, $L_*/K_*$ is spherical. For the second space we 
have to check condition (iii) of Theorem 1. 
We have $K_*(\gt m)=(\U_1)^n$. 
As a $K_*(\gt m)$-module $\gt
n=\gt v_1\oplus...\oplus\gt v_n\oplus\rl$, where $\gt v_i=\rl^2$
for every $1\le i\le n$. 
Each $\gt v_i$ is acted upon by its own $\U_1$.
Note that $[\gt v_i, \gt v_j]=0$ for $i\ne j$. For $K_*(\gt
m)$-invariants we have $S(\gt
n)^{K_*(\gt m)}=\rl[t_1,...t_n,\xi]$, where $t_i$ is the quadratic
$\U_1$ invariant in $S^2(\gt v_i)$ and $\xi\in\gt n'$.
Evidently, $t_i$ and  $t_j$
commute as elements of the  Poisson algebra $S(\gt n)$, and $\xi$
lies in the centre of $S(\gt n)$. 
\end{ex}

Note that the homogeneous space 
$((\cp^n\leftthreetimes\U_n)\times\SU_n)/\U_n$ is a central 
reduction of 
$((H_n\leftthreetimes\U_n)\times\SU_n)/\U_n$ and it is not maximal.   

\begin{thm} Suppose $X=(N\leftthreetimes L)/K$
is a maximal principal indecomposable commutative homogeneous space. 
Then either $X$ is one of the spaces 
$((\rl^n\leftthreetimes\SO_n)\times\SO_n)/\SO_n$,
$((H_n\leftthreetimes\U_n)\times\SU_n)/\U_n$ 
or each
non-commutative simple normal subgroup $K_1\ne{\rm SU}_2$
of $K$ is contained in either $P$ or $L^\diamond$. 
\end{thm}

\begin{proof} Essentially, this is a corollary of Lemma~\ref{9}.
Let $K_1\ne{\rm SU}_2$ be a non-commutative simple normal subgroup
of $K$ that
is not contained in either $L^\diamond$ or $P$. Then 
either $K_1=\SO_n$ or $K_1=\SU_n$ and there are
$P_1\cong L_1^\diamond\cong K_1$ such that 
$K_1\subset P_1\times L_1^\diamond$. Besides, 
$\pi_1^K(K_*)=\SO_{n-1}$ or $\pi_1^K(K_*)=\U_{n-1}$, depending on $K_1$.
According to \cite{al},
$L_1^\diamond$ can act non-trivially only on its simplest module
$V$ ($V=\rl^n$ or $V=\cp^n$) 
and $V\subset\gt n$ is an $L$-invariant subspace.
Set $C_V=Z(L)\cap\GL(V)$. Because $X$ is principal, 
$C_V\subset K$ and if $C_V$ is trivial, 
then $\pi_1^K(K_*)=(K_1)_*(V)$.   
For $K_1^\diamond=\SU_n$ we have 
$(\SU_n)_*(\cp^n)=\SU_{n-1}$, so $C_V=\U_1$, 
$(\U_n)_*(\cp^n)\cong\pi_1^K(K_*)=\U_{n-1}$. 
In case $K_1=\SO_n$
the group $C_V$ is trivial. 

Denote by $\gt n_1:=V+[V,V]$ the Lie subalgebra generated
by $V$, by $N_1\subset N$ the corresponding 
connected subgroup. 
If $K_1=\SO_n$,
then $(K_1)_*(\gt m)=(\U_1)^{[n/2]}$. 
According to condition (iii) of Theorem 1, 
$\gt n_1$ has to be commutative, i.e, $\gt n_1=V$. 
If $K_1=\SU_n$, 
then $\gt n_1$ can be either a Heisenberg or a commutative 
algebra. 

Assume that $X\ne((N_1\leftthreetimes(L_1\times C_V))\times P_1)/
(K_1\times C_V)$. 
Then $L=(L_1\times C_V\times P_1)\times F$, 
$K=(K_1\times C_V)\times H$, 
where $H\subset F$.
Let $\gt a$ be 
an $L$ invariant complement of $\gt n_1$ in $\gt n$.
Recall that the actions $(L_1\times C_V\times P_1):\gt a$ 
and  $F:\gt n_1$ are trivial. 

The remaining part (end) of the proof is the same as in Theorem 3.
If $\gt a$ is a subalgebra (ideal), then $X$ is decomposable. 
Assume that $[\gt a,\gt a]\not\subset\gt a$. The action 
$L_1:\gt a$ is trivial, hence $[\gt a, \gt a]\cap V=0$.    
Thus $[V,V]\subset [\gt a,\gt a]$ and $X$ is not maximal.
\end{proof}

Let $K_1=\SU_2$ be a normal subgroup of $K$.
Suppose it has a non-trivial projections onto
$P_1$ and $L_1^\diamond$. If $L_1^\diamond\ne\SU_2$,
then, as we can see from Table 2a, $L_1^\diamond\cong\Spin_8$.
But as was already mentioned, the pair $(\SU_2\times\Spin_8,
\SU_2\times\Sp_2)$ is not spherical. So $L_1^\diamond=\SU_2$.

If $\pi_1^K(K_*)\ne K_1$, i.e.,  $\pi_1^K(K_*)^0=\U_1$,
then $K_1\subset P_1\times L_1^\diamond$ and $P_1=\SU_2$.
But if  $\pi_1^K(K_*)=K_1$ (and this can be the case),
then $P_1$ can be larger and $K_1$ can have a non-trivial projection onto
some other simple factor
$P_2$ or $L^\diamond_2=\SU_2$.

\begin{ex} Let $\Sp_{m-1,1}$ be a non-compact 
real form of $\Sp_{2m}(\cp)$.
Set $P:=\Sp_{m-1,1}\times\Sp_l$,
$L^\diamond:=\Sp_1\times\Sp_n$,
$K:=\Sp_{m-1}\times\Sp_{l-1}\times\Sp_1\times\Sp_n$ and take for
$N$ a commutative group $\qv^n$.
The inclusions and actions are illustrated by the following diagram.
 \begin{center}
\begin{tabular}{c}
 \xymatrix@R-3mm@C-1mm{
                  {\Sp_{m-1,1}} & {\Sp_l}  & {\Sp_1} \ar[drrr] & {\Sp_n} \ar[drr]
   & &\\
                  {Sp_{m-1}} \ar@{-}[u] & \Sp_{l-1} \ar@{-}[u] & \Sp_1
      \ar@{-}[u] \ar@{-}[ul] \ar@{-}[ull] &
                          {\Sp_n} \ar@{-}[u] & &  {\qv^n} \\ } \\
\end{tabular}
\end{center}

\noindent
The homogeneous space $((N\leftthreetimes L^\diamond)\times P)/K$
is commutative.
Here $L_*=\Sp_{m-1,1}\times\Sp_l\times\Sp_1\times\Sp_{n-1}$ and
$K_*=\Sp_{m-1}\times\Sp_{l-1}\times\Sp_1\times\Sp_{n-1}$.
\end{ex}

To avoid complicated technical details concerning actions
and inclusions of normal subgroup isomorphic to $\Sp_1$ we
impose on $X$ a condition of $\Sp_1$-saturation.

\par\medskip

\begin{center}
 {\large \bf 5.} $\Sp_1${\large \bf-saturated spaces }
\end{center}

Let $X=(N\leftthreetimes L)/K$ be a commutative homogeneous
space.  
Let $L_i$ be a simple direct factor of $L$. By our assumptions 
$L$ is a product 
$L=Z(L)\times L_i\times L^i$, where $L^i$ contains all 
direct factor $L_j$ with $j\ne i$.   

\begin{df} A commutative homogeneous space $X$ is called  {\it
$\Sp_1$-saturated}, if

(1) any normal subgroup $K_1\cong\SU_2$ of $K$ is contained
in either $P$ or $L^\diamond$;

(2) if a simple direct factor $L_i$ 
is not contained in $P$ and $\pi_i(L_*)=L_i$, then
$L_i\subset K$;

(3) if there is an $L$-invariant subspace
$\gt w_j\subset(\gt n/\gt n')$ such that
for some $L_i$ the action $L_i:\gt w_j$ is non-trivial and
the action $Z(L)\times L^i:\gt w_j$ is irreducible,
then $L_i$ acts on
$(\gt n/\gt n')/\gt w_j$ trivially.
\end{df}

\begin{ex}\label{derevo}
Suppose that we have a linear action of a connected compact 
group $F=\Sp_1\times \check F$ on a vector space $V$ and
$F=\check F F_*(V)$, i.e., $\rl[V]^F=\rl[V]^{\check F}$. 
Then we can construct several 
non-$\Sp_1$-saturated commutative 
homogeneous spaces, for instance, 
$(V\leftthreetimes F)/(\U_1\times\check F)$, 
$((V\leftthreetimes F)\times\Sp_1)/(\check F\times\Sp_1)$, 
$((V\leftthreetimes F)\times\Sp_m)/(\check F\times\Sp_1\times\Sp_{m-1})$,
$((V\leftthreetimes F)\times\Sp_{m,1})/(\check F\times\Sp_1\times\Sp_m)$,
where 
$V$ is regarded as a simply connected abelian group.  

Consider a rooted tree with vertices $0,1,\dots, q$, where $0$ is the root.
To each vertex $i$ we attach a positive integer $d(i)$. 
Assume that $d(0)=1$. 
Let $F$ be a product of $\Sp_{d(i)}$ over all vertices, 
and $\check F$ be a product of $\Sp_{d(i)}$ over all vertices 
except the root. 
To each edge $(i,j)$ we attach the vector space
$\qv^{d(i)}\otimes_\qv\qv^{d(j)}$.
Let $V$ be a direct sum of all these spaces.
The group $\Sp_{d(i)}$ naturally acts 
on the first factor in $\qv^{d(i)}\otimes(\bigoplus\qv^{d(j)})$, 
where  the sum is taken over all $j$ connected with $i$. 
For example, a tree with two vertices 
corresponds to a linear representation $\Sp_1\times\Sp_{d(1)}:\qv^{d(1)}$.    
We can calculate $F_*(V)$ consecutively,
descending each time one level in the tree and verifying that 
$F=\check F F_*(V)$.  
Using Lemma~\ref{4} and some basic facts 
concerning representations of symplectic algebras one can prove that
each triple $(F, \check F,  V)$ with 
$\rl[V]^F=\rl[V]^{\check F}$ corresponds to a tree described above.
\end{ex}

This description is complicated, but the commutative 
spaces obtained do not differ much from either reducible ones or 
spaces of Euclidian type. 
Example~\ref{derevo} is just the beginning of another
long story, which will be  
considered elsewhere. 

\begin{ex}\label{sp1} Set 
$X=((N\leftthreetimes(\Sp_n\times\Sp_1))\times\Sp_1)/
(\Sp_n\times\Sp_1)$, where $\gt n=\mathbb H^n\oplus\mathbb H_0$ 
is a two-step nilpotent non-commutative Lie algebra with 
$[\qv^n,\qv^n]=\qv_0$, $\mathbb H_0$
is the space of purely imaginary quaternions,
the normal subgroup $\Sp_1$ of $K$ is the diagonal of the product
$\Sp_1\times\Sp_1$. Here $\qv^n=\qv^n\otimes_{\qv}\qv$, 
where $\Sp_n$ acts on $\qv^n$ and $\Sp_1$ acts on $\qv^1$; 
$\qv_0\cong\gt{sp}_1$ as an $L$-module, i.e., $\Sp_n$ acts
on it trivially and $\Sp_1$ via adjoint representation.   

Evidently, $X=(N\leftthreetimes L)/K$ is not
$\Sp_1$-saturated. We show that it is commutative.
First we compute the generic stabiliser $L_*$. 
Recall that $(\Sp_n\times\Sp_1)_*(\qv^n)=
\Sp_{n-1}\times\Sp_1$.
Clearly $L_*(\qv_0)=\Sp_n\times\U_1\times\Sp_1$,
$(\Sp_n\times\U_1)_*(\qv^n)=\Sp_{n-1}\times\U_1$. 
We have 
$L_*=\Sp_{n-1}\times\U_1\times\Sp_1$, 
$K_*=K\cap L_*=\Sp_{n-1}\times\U_1$, 
$K_*(\gt m)=\Sp_n\times((\Sp_1)_*(\gt{sp}_1))=
\Sp_n\times\U_1$.
One can easily verify conditions (i) and (ii) of Theorem 1. 
Tables of \cite{Vin} and \cite{Vin2}
shows that (iii) is also satisfied.
\end{ex}


Let $X$ be a non $\Sp_1$-saturated commutative homogeneous space.
It can be made $\Sp_1$-saturated by enlarging $K$, $L$ and possibly
$N$ too.
For instance,
if a  simple factor $\Sp_1$ of $K$ has non-trivial
projections onto $P$ and $L^\diamond$,
then we replace $P$ by $P\times\Sp_1$ or
$P\times\Sp_1\times\Sp_1$ (the second replacement is
needed if $\Sp_1$ has non-trivial projections onto
two simple factors of $P$). The group
$K$ is replaced by $K\times\Sp_1$. Starting with the
commutative spaces from Example~\ref{sp1}
we construct an $\Sp_1$-saturated commutative homogeneous space
 $(\Sp_1\times\Sp_1/\Sp_1)\times(N\leftthreetimes K/K)$, where
 $K=\Sp_1\times\Sp_n$ and  $N$ is the same as before.

\begin{ex}
Set $L=K=\Sp_n\times\Sp_1\times\Sp_m$,
$\gt n=\qv^n\oplus\qv^m\oplus\qv_0$, where both algebras
$\qv^n$ and $\qv^m$ are not commutative and
$[\qv^n, \qv^n]=[\qv^m, \qv^m]=\qv_0$.
We have
$S(\gt g/\gt k)^K=S(\gt n)^K=\rl[\xi_1, \xi_2, \eta]$,
where $\xi_1\in S^2(\qv^n)^{\Sp_n}$,
$\xi_2\in S^2(\qv^m)^{\Sp_m}$,
$\eta\in S^2(\qv_0)^{\Sp_1}$, so
the corresponding homogeneous space
$(N\leftthreetimes K)/K$ is commutative.
The third condition of Definition 6 is not fulfilled.
If we want to enlarge $L$, we also need to enlarge
$N$. As a $\Sp_1$-saturation we have a product of two
commutative spaces
$(N_i\leftthreetimes K_i)/K_i$,
where $\gt n_1=\qv^n\oplus\qv_0$,
$K_1=\Sp_n\times\Sp_1$;
$\gt n_2=\qv^m\oplus\qv_0$,
$K_2=\Sp_m\times\Sp_1$.
\end{ex}


The procedure that is inverse  to  $\Sp_1$-saturation 
can have steps of three different types. 
First, one simple factor $\Sp_1$ of $K$ is replaced by
$\U_1$; second, two of three simple factors $\Sp_1$ of $K$ are replaced by
the diagonal of their product; third, several simple factors $\Sp_1$ of
$L$
are replaced by
the diagonal of their product, $K$ is replaced by the intersection
of the former $K$ and new $L$ and probably $N$ is decreased.

According to Lemma \ref{8}, if $(L, K)$ contains any of the
``$\SU_2$'' pairs condition (2) of Definition 6 is not
satisfied. So an  $\Sp_1$-saturated maximal principal commutative
homogeneous space is a product of the spaces from
Table 2b, spaces of the type
$((\rl^n\leftthreetimes\SO_n)\times\SO_n)/\SO_n$,
$((H_n\leftthreetimes\U_n)\times\SU_n)/\U_n$, spherical
homogeneous space of a reductive Lie group and a space of
Heisenberg type. 

\par\medskip

\begin{center}
{\large \bf 6. Commutative homogeneous spaces of Heisenberg type}
\end{center}

In this section we consider homogeneous spaces of the type
$(K\rightthreetimes N)/K$.
In this case $S(\gt g/\gt k)^K=S(\gt n)^K$.
We assume that $\gt n$ is not commutative.

Let us decompose $\gt n/\gt n'$ into a sum of irreducible
$K$-invariant subspaces,
namely $\gt n/\gt n'=\gt w_1\oplus ... \oplus\gt w_p$. As was proved
in \cite{Vin},
if $X$ is commutative then $[\gt w_i, \gt w_j]=0$ for $i\ne j$,
also $[\gt w_i, \gt w_i]=0$ if there is $j\ne i$
such that $\gt w_i\cong\gt w_j$ as a $K$-module. Denote by
$\gt n_i:=\gt w_i\oplus [\gt w_i, \gt w_i]$ the subalgebra generated
by $\gt w_i$. Let $\gt v^i$ stand for a $K$-invariant complement of
$\gt n_i$ in
$\gt n$ and set $K^i:=K_*(\gt v^i)$.

\begin{thm}\label{heiz}(\cite[Theorem 1]{Acta})
In the above notation, $G/K$ is commutative \iff
each Poisson algebra $S(\gt n_i)^{K^i}$ is commutative.
\end{thm}

Note that the irreducibility of $\gt w_i$ is not important here. The statement
of the theorem remains true
for any $K$-invariant subspace $\gt w\subset\gt n/\gt n'$.

For convenience of the reader we present here the
classification results of \cite{Vin} and
\cite{Vin2}. All maximal commutative homogeneous
spaces of Heisenberg type with $\gt n/\gt n'$ being an
irreducible $K$-module and $\dim\gt n'>1$ 
are listed in Table 3.
The following notation is used:

$\gt n=\gt w\oplus\gt z$,
where $\gt z=\gt n'$ is the centre of $\gt n$;


$\mathbb H_0$ is the space of purely imaginary quaternions;

$\cp^m\otimes\qw^n$ is the tensor product over $\cp$;

$\qw^m\otimes\qw^n$ is the tensor product over $\qw$;

$H\Lambda^2\mathbb D^n$, where $\mathbb D=\cp$ or $\qw$,
is the skew-Hermitian square of $\mathbb D$;

$HS_0^2\mathbb H^n$ is the space of Hermitian quaternion matrices
of order $n$ with zero trace.

\noindent
For all cases in Table 3 the commutation operation $\gt w\times\gt w\mapsto\gt
z$
is uniquely determined by the condition of $K$-equivariance.
Notation $(\U_1\cdot) F$ means
that $K$ can be either $F$ or $\U_1\cdot F$. Cases in which $\U_1$ is necessary
are indicated in the column
``$\U_1$''. Some spaces are not always maximal. 
This is indicated in the column
``max''.

\par\smallskip

\begin{tabular}{|c|c|c|c|c|c|}
\multicolumn{6}{c}{Table 3.}\\
  \hline
  & $K$ & $\gt w$ & $\gt z$ & $\U_1$ & max\\
\hline
1 & $\SO_n$ & $\rl^n$ & $\Lambda^2\rl^n=\gt{so}_n$ &&\\
2 & $\Spin_7$ & $\rl^8$ & $\rl^7$ &&\\
3 & $G_2$ & $\rl^7$ & $\rl^7$ &&\\
\hline
4 & $(\U_1\cdot) \SU_n$ & $\cp^n$ & $\Lambda^2\cp^n \oplus\rl$ &&\\
&($n$ even)&&&&\\
5 & $(\U_1\cdot) \SU_n$ & $\cp^n$ & $\Lambda^2\cp^n$ &&\\
&($n$ odd)&&&&\\
6 & $\U_n$ & $\cp^n$ & $H\Lambda^2\cp^n = \gt u_n$ &&\\
7 & $(\U_1\cdot) \Sp_n$ & $\qv^n$ & $HS^2_0\qv^n\oplus\qv_0$ &&\\
8 & $\U_1\cdot\Spin_7$ & $\cp^8$ & $\rl^7\oplus\rl $ &&\\
\hline
9 & $\Sp_1\times\Sp_n$ & $\qv^n$ & $\qv_0=\gt{sp}_1$ && $n\geqslant 2$\\
10 & $\Sp_2\times\Sp_n$ & $\qv^2\otimes\qv^n$ & $H\Lambda^2\qv^2=\gt{sp}_2$
&&\\
\hline
11 & $(\U_1\cdot)\SU_2\times \SU_n $ & $\cp^2\otimes\cp^n $ &
$H\Lambda^2\cp^2=\gt u_2$
    & $n=2$ &\\
12 & $\U_2 \times \Sp_n $ & $\cp^2\otimes\qv^n $ & $ H\Lambda^2\cp^2=\gt u_2$
&&\\
\hline
\end{tabular}
\par\smallskip

In the general case, the classification 
of maximal principal $\Sp_1$-saturated 
commutative spaces of Heisenberg type
is being done in the following way.
If $(N\leftthreetimes K)/K$ is a commutative homogeneous space of
non-Euclidian type, then there is a non-commutative subspace
$\gt w_1\subset({\gt n}/{\gt n}')$.
Denote by $K_e$ the maximal
connected subgroup of $K$ acting on $\gt w_1$ effectively
and by $\pi_e$ the projection onto $K_e$ in $K$.
Recall that $\gt n_1:=\gt w_1\oplus[\gt w_1,\gt w_1]$. 

If $\gt n_1$ is not a Heisenberg algebra (i.e., $\dim\gt n_1'>1$), 
then the pair $(K_e, \gt n_1)$ is 
a central reduction of some pair from Table 3. 
If $\gt n_1$ is a Heisenberg algebra, then $(K_e, \gt n_1)$
corresponds
to an irreducible spherical representation $K_e(\cp):W$, 
in a sense that $\gt w_1$ is a $K_e$-invariant 
real form of $W\oplus W^*$. According to Kac's list
\cite{Kac}, 
there are 14 such cases. 

For any commutative homogeneous space $(N_1\leftthreetimes K_e)/K_e$,
where $\gt n_1/\gt n_1'=\gt w_1$, we have to find out
if it arise in the aforementioned way from some larger
commutative homogeneous space $(\tilde N\leftthreetimes K)/K$,
and if so, list all of them.
Note that $K^1$ acts on $\gt n_1$ as $\pi_e(K^1)$.
Due to lemma \ref{4} and the third condition
of Definition 6, $\pi_e(K^1)$ is a proper subgroup of
$K_e$.

As was proved in \cite{be-ra}, homogeneous space
$(N\leftthreetimes F)/F$, where $\gt
n=\rl^n\oplus\Lambda^2\rl^n$, is not commutative for any
proper subgroup $F\subset\SO_n$. So if $\gt n_1$ corresponds
to the first row of Table 3, then, by Lemma \ref{4}, $\gt
n=\gt n_1$. 

We say that the action $K:\gt n$ is {\it commutative}
if the corresponding homogeneous space
$(N\leftthreetimes K)/K$ is commutative.

\begin{lm}\label{nepr}
Suppose that $K_e=K_e'\times\Upp_1$,
$[\gt w_1, \gt w_1]\ne 0$ and $\gt w_1=W\otimes\rl^2$, where
$K_e'$ acts on $W$ and $\Upp_1$ on $\rl^2$. Let $F$ be a
proper subgroup of $K_e'$. If the action $F:W$ is reducible
then $(N_1\leftthreetimes (F\times\Upp_1))/(F\times\Upp_1)$
is not commutative.
\end{lm}

\begin{proof} We show that the action
of $H=(\SO_n\times\SO_m)\times\SO_2$ on $\gt n_1
\cong(\rl^n\oplus\rl^m)\otimes\rl^2\oplus[\gt w_1, \gt w_1]$
cannot be commutative. Assume that it is commutative and
apply Theorem 5. We have
$H_*(\rl^n\otimes\rl^2)^0=\SO_{n-1}\times\SO_m$. The
subspace $\rl^m\otimes\rl^2$ is a sum of two isomorphic
$\SO_{n-1}\times\SO_m$-modules. Hence, $\rl^m\otimes\rl^2$
is a commutative subalgebra of $\gt n_1$. This can happen
only if $[\gt w_1, \gt w_1]=0$.
\end{proof}


\begin{lm}\label{neprt3}
Let $(N\leftthreetimes K)/K$ be a commutative homogeneous space
from Table 3, but not from the second, third or ninth row.
Suppose a subgroup $F\subset K$ acts on $\gt n/\gt n'$
reducibly. Then $(N\leftthreetimes F)/F$ is not commutative.
\end{lm}

\begin{proof} Assume that $(N\leftthreetimes F)/F$
is commutative. Then due to Proposition 15 of \cite{Vin}
there are at list two subspaces $V_1, V_2\subset\gt n/\gt n'$,
such that $V_1\oplus V_2=\gt n/\gt n'$ and
$[V_1, V_2]=0$. Evidently, this is not true in cases 1, 4, 5, 6.
For the same reason, in cases 10, 11 and 12
$F$ contains the first simple factor of $K$,
either $\Sp_2$ or $\SU_2$.

In the seventh case $F$ has to be a subgroup
of $Sp_m\times\Sp_{n-m}$. Subspaces $\qw^m$ and $\qw^{n-m}$
do not commute with each other.
For the eighth case we apply Lemma \ref{nepr}.

In case 10, we have $F\subset\Sp_2\times\Sp_m\times\Sp_{n-m}$,
$F_*(\qw^2\otimes\qw^m)\subset\Sp_1\times\Sp_1\times\Sp_n$.
The subspace
$\qw^2\otimes\qw^{n-m}$ is a sum of two isomorphic
$F_*(\qw^2\otimes\qw^m)$-modules. Hence,
$\qw^2\otimes\qw^{n-m}$ is a commutative subalgebra of
$\gt n$. But this is not the case.

In case 11, $F$ is a subgroup of either $\SU_2\times\U_m\times\U_{n-m}$
or $\SU_2\times\Sp_{m/2}$ for even $m$.
The first case is just the same as the previous one.
For the second case note that as a $\SU_2\times\Sp_{m/2}$-module
$\cp^2\otimes\cp^n\cong\qw^{m/2}\oplus\qw^{m/2}$.
The 12-th case is exactly the same.
\end{proof}


{\cled }(of the proof). Let
$(N\leftthreetimes K)/K$ be a non-Euclidian central reduction
of a homogeneous space from row 8, 11 or 12 of Table 3.
Then $(N\leftthreetimes F)/F$ is not commutative for
any proper subgroup of $F\subset K$ acting on $\gt n/\gt n'$
reducibly.

Note that this statement is not true for
a central reduction
$(\U_n, \cp^n\oplus\rl)$ of the pair from the 6-th row of Table 3. 

\begin{lm}\label{tor}
Suppose $\pi_e(K^1)=(\Upp_1)^n$ and a homogeneous space
$(N_1\leftthreetimes K^1)/K^1$ is commutative.
Then $\gt n_1=\rl^{2n}\oplus\rl$, $K_e=\Upp_n$.
\end{lm}

\begin{proof} An irreducible representation of
$\Upp_1$ on a real vector space is either trivial ($\rl$) or
$\rl^2$. If  $\gt w_1$ is the sum of more than
$n$ $K^1$-invariant summands, then two of them
are isomorphic and there is a non-zero $\eta\in\gt w_1$ such that
$[\eta, \gt w_1]=0$.
But $K_e\eta=\gt w_1\subset\gt z({\gt n}_1)$.
By the same reason $\gt w_1^{(\U_1)^n}=0$.
Because the action $(U_1)^n:\gt w_1$ is locally
effective, $\gt w_1=\rl^{2n}$.
We have $\Lambda^2\rl^2=\rl$, hence
$K^1$ acts on ${\gt n}_1'$ trivially.
Each element of
$K_e$ is contained in some maximal torus,
that is up to conjugation in $\pi_e(K^1)$.
Hence $K_e$ acts on ${\gt n}_1'$ trivially and
$K_e\subset\Upp_n$. The group $\Upp_{n}$ has no proper subgroups
of rank $n$ acting on
$\rl^{2n}$ irreducibly. Thus we have
$K_e=\Upp_n$, $\gt n_1'=(\Lambda^2\rl^n)^{\Upp_n}\cong\rl$.
\end{proof}

From now on, let $(N\leftthreetimes K)/K$ be an indecomposable
maximal $\Sp_1$-saturated principal commutative space with
$\gt n_1\ne\gt n$. In particular, the connected centre
$Z(K_e)$ of $K_e$ acts on $\gt v^1$ trivially. Let $\gt a\subset\gt n$ 
be a $K$-invariant subalgebra. Clearly, if the action 
$K:\gt n$ is commutative, then $K:\gt a$ is also commutative. 
We assume that $K:\gt n$ is not a ``subaction'' of some larger 
commutative action.     

Assume for the time
being that $K_e'$ is simple and denote it by $K_1$.
Decompose $\gt n$ into a sum of $K$-invariant subspaces
${\gt n}={\gt n}_1\oplus V_2\otimes_{\mathbb D_2}
V^2\oplus... \oplus V_q\otimes_{\mathbb D_q} V^q\oplus
V_{\rm tr}$, where $V_i$ are pairwise non-isomorphic
irreducible non-trivial $K_1$-modules, $V_{\rm tr}$ and
$V^i$ are trivial $K_1$-modules and all the other simple
normal subgroups of $K$ act on $V_i$ trivially. First we
have to describe possible $V_i$, then dimensions of $V^i$,
afterwards the actions $K:\bigoplus\limits_{i=2}^{q}
V_i\otimes V^i$ and $K:V_{\rm tr}$. Once again we use
\'Elashvili's classification \cite{al}. Lemma \ref{tor}
tells us that the adjoint representation can be among $V_i$
only for $(K_e, \gt n_1)=(\Upp_n, \cp^n\oplus\rl)$. Note
that if $V_{\rm tr}$ is not contained in $\gt n'$, then
there is a simple factor of $K$ acting non-trivially both on
$\bigoplus\limits_{i=2}^{q} V_i\otimes V^i$ and $V_{\rm
tr}$.

Recall that $\pi_e(K^1)$ is contained in the product
$Z(K_e)\cdot N(\xi)$, where $N(\xi):=\{k\in K_1|
k\xi\in\mathbb D_i\xi\}$ is the ``normaliser'' of a generic
vector $\xi\in V_i$. Suppose $(K_1)_*(V_i)$ is trivial or
finite. Then $\pi_e(K^1)$ is finite for $\mathbb D_i=\rl$
and commutative for $\mathbb D_i=\cp$. In the following
lemma we prove that for $\mathbb D_i=\qw$ the projection
$\pi_e(K^1)$ is also commutative.

\begin{lm}\label{lsp}
Let $F\subset\Sp_n$, $n\ge 2$ and $(F_*(\mathbb H^n))^0=E$.
Then the image of the generic stabiliser
$(F\times\Sp_m)_*(\mathbb H^n\otimes\mathbb H^m)$ under the
projection on $F$ is commutative.
\end{lm}

\begin{proof} Assume that this is not the case, i.e.,
the image contains $\Sp_1$.
Then this is true not only for generic
points, but for all of them, in particular, for decomposable
vectors. In other words, the restriction
$F|_{\xi\qw}$ for $\xi\ne 0, \xi\in\qw^n$ contains $\Sp_1$.
Due to Lemma \ref{4}, we have $F=\Sp_n$. But then
$F_*(\mathbb H^n)=\Sp_{n-1}$.
\end{proof}

\begin{ex} Let $(K_e, \gt n_1)$ be a pair from the second
row of Table 3. We show that
$\gt n\subset\gt n_1\oplus\rl^7\otimes\rl^2$.
All representations of $\Spin_7$ are orthogonal, so
here all $\mathbb D_i$ equal $\rl$. The group $\Spin_7$
has only three irreducible representations with infinite
generic stabiliser, namely $\gt{so_7}$, $\rl^7$ and $\rl^8$. If
$V_i=\rl^8$ for some $i$, then $K^1$ has a non-zero invariant in
$\gt w_1$, which commute with $\gt w_1$. This is a
contradiction.
Thus $\gt n=\gt n_1\oplus\rl^7\otimes V^2\oplus V_{\rm tr}$.
If $\dim V^2\ge 3$, then $\pi_e(K^1)\subset\SU_2\times\SU_2$.
But the action $\SU_2\times\SU_2:(\cp^2\oplus\cp^2)\oplus\mathbb R^7$
is not commutative. In case $\dim V^2=2$ we have
$\pi_e(K^1)=\Spin_5=\Sp_2$. The pair
$(\Sp_2, \qw^2\oplus\rl^7)$ is a central reduction
of the pair from the 7-th row of Table 3 with $n=2$ 
by a subgroup corresponding to $\qv_0$
(here $HS_0^2\qv^2\cong\rl^7$ as a $\Sp_2$-module). 

We have seen that $\dim V^2\le 2$, so no
simple normal subgroup of
$K$ except $\Spin_7$ acts non-trivially on $V_2\otimes V^2$.
Assume that there is a non-commutative subspace
$\gt w_2\subset\gt n$.
It is either $\rl^7$ or $\rl^7\otimes\rl^2$.
The first case is not possible, because
$\Lambda^2\rl^7\cong\gt{so}_7$ as a $\Spin_7$-module.
In the second case we apply
Theorem 5 to $\gt w_2=\rl^7\otimes\rl^2$.
Recall that, $K=\Spin_7$ or
$K=\Spin_7\times\SO_2$,
$(\Spin_7)_*(\rl^8\oplus\rl^7)\subset\Spin_6$.
By Lemma 11, the action $\SO_6\times\SO_2:
\gt w_2\oplus[\gt w_2, \gt w_2]$ is commutative only if
$[\gt w_2, \gt w_2]=0$.
The corresponding commutative space is indicated in the
13-th row of Table 4.
\end{ex}

There are 9
cases in which $K_e$ has two simple factors.
Namely, 4 last rows of
Table 3 and there central reductions, 
and 3 irreducible spherical representations:
$(\cp^*\cdot)\SL_m\times\SL_n:\cp^m{\otimes}\cp^n$,
$(\cp^*\cdot)\SL_n\times\Sp_4:\cp^n{\otimes}\cp^4$,
$\GL_3\times\Sp_n:\cp^3{\otimes}\cp^{2n}$. We have to
carry the same procedure for both simple normal subgroups of $K$.

For convenience of the reader, we list here all irreducible
representations of $\gt{su}_n$ with non-trivial generic
stabiliser. They are described by the highest weights of the
complexifications $V(\cp)$ with respect to $\gt{sl}_n$.

\begin{center}
Table $A_{n-1}$.

 \nopagebreak
\begin{tabular}{|>{$}c<{$}|>{$}c<{$}|>{$}c<{$}|}
\hline
n & \mbox{representation } & (\SU_n)_*(V)^0 \\
\hline
 & R(\varpi_1)\oplus R(\varpi_1)^* & \SU_{n-1} \\
\hline
 & R(\varpi_2)\oplus R(\varpi_2)^* & (\SU_2)^{[n/2]} \\
\hline
 & R(2\varpi_1)\oplus R(2\varpi_1)^* & (\Upp_1)^{[n/2]} \\
\hline
 & R(\varpi_1+\varpi_1^*) & (\Upp_1)^{n-1} \\
\hline
 4 & R(\varpi_2) & \Sp_2 \\
\hline
 6 & 2R(\varpi_3) & (\Upp_1)^2 \\
\hline
\end{tabular}
\end{center}

\noindent
Note that the action
$(\SU_n)_*(V):\cp^n$ is irreducible
only in one case: $n=4$, $V=R(\varphi_2)\cong\rl^6$.

For almost all pairs from Table 3, the action
$K^1:\gt w_1$ have to be irreducible. This leaves only a few possibilities
for $V_i$. The obtained commutative spaces 
are listed in rows 2, 4, 5 of Table 4. 

The Lie group $G_2$ has only 
two irreducible representations with non-trivial generic stabiliser,
namely adjoint one and $\rl^7$.
Thus, if $(K_e, \gt n_1)$ is the pair from the 3-d row of Table 3, 
then $K=K_e$ and $\gt n=\gt n_1$. 

Calculations in cases 
$((\U_1)\cdot\Sp_n, \qv^n\oplus\qv_0)$,
$((\U_1)\cdot\Sp_n, \qv^n\oplus\rl)$
and $(\Sp_1\cdot\Sp_n, \qv^n\oplus\gt{sp}_1)$
do not differ much. By our assumptions 
subgroups $\U_1$ and $\Sp_1$ act on $\gt v^1$ trivially.
The result is given in rows 9, 10, 11 and 12 of Table 4. 
  
If $\gt n'$ is a trivial $K$-module, the calculations are even simpler.
However, we have more such cases. Assume that
$\gt n_1$ is not commutative. Recall that
$\gt w_1(\cp)=W_1\oplus W_1^*$ as a $K$-module.
We have to check whether the action $\pi_e(K^1): W_1$
is spherical or not. We will consider one example in full
details.

\begin{lm}\label{su2}
Let $(N\leftthreetimes K)/K$, where $K=\SU_2\times F$ be a
principal $\Sp_1$-saturated commutative space. Suppose that
the image of the projection of $K_*({\gt n})$ on $\SU_2$
contains $\Upp_1$. Then there are three possibilities:
$K=\SU_2$, ${\gt n}=\rl^3$; $F=(\Spp)U_s$ or
$F=(\Upp_1\cdot)\Sp_{s/2}$, $\gt
n=\cp^2\otimes\cp^s\oplus\rl$; $F=(\Spp)U_4$, $\gt
n=(\cp^2\otimes\cp^4\oplus\rl)\oplus \rl^6$.
\end{lm}

\begin{proof} Note that the intersection of
two distinct general subgroups
$\Upp_1\subset\SU_2$ is finite. Hence,
$\gt n$ contains only one irreducible $K$-invariant
subspace on which $\SU_2$
acts non-trivially. The Lie algebra $\gt{su}_2$
has only two non-trivial representations
with a non-trivial generic  normaliser, namely
$\rl^3$ and $\cp^2$. The first one is orthogonal and can not
be contained in $\gt n$ in the form
$\rl^3\otimes\rl^s$ for $s>1$.
Assume that $\cp^2\otimes\cp^s\subset{\gt n}$.  The group
$F$ acts on a generic subspace $\cp^2\subset\cp^s$
as $\Upp_1$ or as  $\SU_2$. It can be easily
verified that the first case is not possible.
In the second one
$F$ and $\Upp_m$ have the same orbits in
$\cp^m$, i.e., $\Upp_m=F\Upp_{m-1}$. According
to the classification \cite{On} , there are only 4
listed in the lemma possibilities for
$F$. Suppose $F$ acts (non-trivially) on some other
$K$-invariant subspace $V\subset\gt n$.
Then $\cp^s$ is an irreducible $F_*(V)$-module.
This can only happen for $s=4$,
$F=(\Spp)U_4$ and  $V=\rl^6$.
\end{proof}

\begin{ex}\label{A}
We describe all possible principal $\Sp_1$-saturated
maximal commutative
pairs $(K, \gt n)$ for $(K_e, \gt n_1)=((\Upp_1\cdot)\SU_n,
\cp^n\oplus\rl)$.

First assume that $n=2$. The action
$\Upp_1\cdot\Upp_1:\cp^2\oplus\rl$ is commutative. Hence we
can "replace" $\SU_2$ by $\Upp_1$. This yields the three 
possibilities from Lemma~\ref{su2}.

Suppose $n>2$. From Table $A_{n-1}$ we know that $\gt{su}_n$
has the following irreducible representations with
non-trivial generic stabiliser: $\gt{su}_n$, $\cp^n$,
$\Lambda^2\cp^n$ for $n>4$ and $\rl^6$ for $n=4$. All of
them can occur in $\gt v^1$ as $V_i$. The first one is
orthogonal, so for $s>1$ the space $\gt{su}_n\otimes\rl^s$
could not be a subspace of $\gt v^1$. Moreover, if
$V_2=\gt{su}_n$, then $\gt v^1=\gt{su}_n$, because a rank of
any proper subalgebra $\gt f\subset\gt{su}_n$ is smaller
then $n-1$. The action $\U_n:\gt n_1\oplus\gt{su_n}$ is
commutative. In the following we assume that $\gt{su}_n$ is
not contained in $\gt v^1$.

Consider the case $n=4$. We have
$\gt v^1=\rl^6\otimes\rl^s\oplus\cp^4\otimes V^3\oplus V_{\rm tr}$.
Note that $(\SU_4)_*(\rl^6\otimes\rl^3)$ is finite, so
$s\le 2$. Also $(\Upp_1\cdot\Upp_4)_*(\rl^6\otimes\rl^2\oplus
\cp^4)=(\Upp_1)^3$ and $\cp^4$ is not a spherical representation
of $(\cp^*)^3$. Hence, for $s=2$ we have
$(K, \gt n)=((\Spp)\Upp_4(\cdot\SO_2),
\cp^4\oplus\rl\oplus\rl^6\otimes\rl^2)$.  Another possible pair
is $(\Upp_4\cdot\Upp_1, (\cp^4\oplus\rl)\oplus\rl^6
\oplus(\cp^4\oplus\rl))$. The rest of this case is
the same for general $n$ and is dealt upon below.

Note that
$(\SU_n)_*(\Lambda^2\cp^n\oplus\Lambda^2\cp^n)=\Upp_1$ and
$(\SU_n)_*(\Lambda^2\cp^n\oplus\cp^n)=\{E\}$. Thus we have
either $\gt n=\gt n_1\oplus(\Lambda^2\cp^n\oplus\rl)$ or
$\gt n=\gt n_1\oplus\cp^n\otimes_{\mathbb D_2} V^2\oplus
V_{\rm tr}$. Here $\mathbb D_2$ equals $\cp$ or $\rl$. If
$\mathbb D_2=\rl$ and $\dim V^2>1$, then $\pi_e(K^1)$ is
contained in
$\U_1\cdot\U_{n-2}\subset\Upp_2\cdot\Upp_{n-2}$. Evidently,
the $\U_1\cdot\U_{n-2}$-module $\cp^n$ is not spherical.
 Hence, $\mathbb D_2=\cp$.

Similar to the proof of Lemma \ref{su2} we obtain that
$K_*(V_{\rm tr})/\SU_n$ acts on $V^2=\cp^s$ as
$(\Spp)\Upp_s$ or $\Sp_{s/2}$ for even $s$. In the second
case (for $s\ge 4$) we have $\pi_e(K^1)\subset
(\Upp_1)^2\cdot\SU_{n-3}$, but the action $(\cp^*)^2:\cp^3$
is not spherical. Thus $V_{\rm tr}=\rl$ is a trivial
$K$-module and $\gt n\subset\gt
n_1\oplus\cp^2\otimes\cp^s\oplus\rl$.
\end{ex}

It is not difficult to show, that 
If $(K_e, \gt n_1)$ is one of the pairs:
$(\U_1\cdot\SO_n, \gt h_n)$ with $n\ne 8$, 
$(\U_n, S^2\cp^n\oplus\rl)$, 
$(\U_1\cdot\Spin_9, \gt h_{16})$, 
$(\U_1\cdot G_2, \gt h_7)$,
$(\U_1\cdot E_26, \gt h_{27})$,
then $K=K_e$, $\gt n=\gt n_1$. One has to use
the list of \cite{al} and for some cases Lemma~\ref{nepr}.  

We will not consider other cases. 
They are similar to Example~\ref{A}.

All commutative homogeneous spaces of Hiesenberg 
type satisfying our restrictions are listed in Table 4. 
 The algebra $\gt n_{\mbox{\small max}}$ is given in the following way.
Each subspace in
parentheses represent a subalgebra $\gt w_i\oplus[\gt w_i, \gt w_i]$.
The spaces given outside parentheses are commutative.
The action of $K$ is uniquely determined by irreducibility and
by Table 3.

Notation $(\SU_n,\Upp_n, \Upp_1\cdot\Sp_{n/2})$
means that this normal subgroup of $K$ can be equal
either of these three groups.
Symbol $\Sp_{n/2}$ has sense only for even $n$.
For $n$ odd, the group
$\Sp_{n/2}$ does not exist, so it cannot be a subgroup of $K$.

\begin{thm} All indecomposable $\Sp_1$-saturated maximal principal
commutative homogeneous spaces $(N\leftthreetimes K)/K$
 with non-commutative $\gt n$ and reducible
$\gt n/\gt n'$ are given in Table 4 in a sense that $\gt n$
is a $K$-invariant subalgebra of $\gt n_{\mbox{\small max}}$ .
\end{thm}

\begin{proof} We have seen above how one can prove that
all such commutative spaces are contained in Table 4. Now we
prove that all listed spaces are commutative. We do it using
Theorem~\ref{heiz} and the list of the spherical representations
given in \cite{kn}. It is proved in \cite{be-ra}, 
\cite{be-ra2} and \cite{L} that
spaces listed in  rows
3,\,7,\,8,\,9,\,15,\,18 and 19 are commutative.

Suppose $\gt n$ contains a commutative $K$-invariant
ideal $\rl^6$. According to Theorem \ref{heiz},
$K:\gt n$ is commutative \iff $K_*(\rl^6):\gt n/\rl^6$
is commutative. For second, 4-th and 5-th rows of Table 4
pairs $K_*(\rl^6):\gt n/\rl^6$ are contained in Table 3, hence,
commutative. For 6-th, 20-th and 23-d rows of Table 4
pairs $K_*(\rl^6):\gt n/\rl^6$ correspond to spherical representations
listed in \cite{be-ra2}, \cite{L}. Analogously, for the 22-d row
pair $K_*(\rl^6\oplus\rl^6):\gt n/(\rl^6\oplus\rl^6)$ corresponds
to a spherical representation.

Let $(N\leftthreetimes L)/K$ be a commutative homogeneous
space. Then the action
$K:\gt n\oplus(\gt l/\gt k)$, where
$\gt l/\gt k$ is a commutative ideal,
is commutative. 
The pairs from the first and 12-th
rows of Table 4 are obtained in such a way from
commutative spaces
$((H_n\leftthreetimes\U_n)\times\SU_n)/\Upp_n$
and $(H_{2n}\leftthreetimes\U_{2n})/\Sp_n$.
In case $[\gt n,\gt n]=0$, we obtain a commutative 
homogeneous space of Euclidian type. Thus essentially these 
are the only non-trivial examples given by this construction. 

In the remaining  eight cases we use Theorem \ref{heiz}.
For instance, take the 11-th row with
$K=\Sp_n\times\Sp_m$. Here
$\gt n$ contains only one non-commutative subspace
$\gt w_1\cong\qv^n$. Set
$d=|n-m|$ and $s=\min(n, m)$. Then
$K^1=K_*(\qv^n\times\qv^m)=(\Sp_1)^d\times\Sp_s$.
Anyway $K^1/(K^1\cap P_1)$
contains $(\Sp_1)^n$. To conclude note that
the action $Sp_1:(\qv\oplus\qv_0)$ is commutative according
to Table 3.
\end{proof}

\begin{center}
\begin{tabular}{|c|>{$}c<{$}|>{$}c<{$}|}
\multicolumn{3}{c}{Table 4.} \\
\hline
& K & {\gt n_{\mbox{\small max}}}\\
\hline
1 & \Upp_n & (\cp^n\oplus\rl)\oplus\gt{su}_n \\
\hline
2 & \Upp_4 & (\cp^4\oplus\Lambda^2\cp^4\oplus\rl)\oplus\rl^6 \\
\hline
3 & \Upp_1\cdot\Upp_n & (\cp^n\oplus\rl)\oplus
(\Lambda^2\cp^n\oplus\rl) \\
\hline
4 & \SU_4 & (\cp^4\oplus HS_0^2\mathbb H^2\oplus\rl)\oplus\rl^6 \\
\hline
5 & \Upp_2\cdot\Upp_4 & (\cp^2\otimes\cp^4\oplus
H\Lambda^2\cp^2)\oplus\rl^6\\
\hline
6 & \SU_4\cdot\Upp_m & (\cp^4\otimes\cp^m\oplus\rl)\oplus\rl^6\\
\hline
7 & \Upp_m\cdot\Upp_n &
(\cp^m\otimes\cp^n\oplus\rl)\oplus(\cp^m\oplus\rl)\\
\hline
8 & \Upp_m\cdot\SU_2\cdot\Upp_p & (\cp^m\otimes\cp^2\oplus\rl)
\oplus(\cp^2\otimes\cp^p\oplus\rl)\\
\hline
9 & \Upp_1\cdot\Upp_1\cdot\Sp_n & (\mathbb H^n\oplus\rl)\oplus
(\mathbb H^n\oplus\rl)\\
\hline
10 & \Sp_n\cdot\Sp_1\cdot(\Sp_1, \Upp_1, \{E\}) & (\mathbb H^n\oplus\mathbb
H_0)\oplus
\mathbb H^1\otimes\mathbb H^n\\
\hline
11 & \Sp_n\cdot(\Sp_1, \Upp_1)\cdot\Sp_m &
(\mathbb H^n\oplus\mathbb H_0)\oplus \mathbb H^n\otimes\mathbb H^m\\
\hline
12 & \Sp_n\cdot(\Sp_1, \Upp_1, \{E\}) & (\mathbb H^n\oplus\mathbb H_0)\oplus
HS_0^2\mathbb H^n\\
\hline
13 & \Spin_7\cdot(\SO_2, \{E\}) & (\rl^8\oplus\rl^7)\oplus \rl^7\otimes\rl^2\\
\hline
14 & \Upp_1\cdot\Spin_7 & (\cp^7\oplus\rl)\oplus \rl^8 \\
\hline
15 & \Upp_1\cdot\Upp_1\cdot\Spin_8 & (\cp^8_+\oplus\rl)\oplus
(\cp^8_-\oplus\rl)\\
\hline
16 & \Upp_1\cdot\Spin_{10} & (\cp^{16}\oplus\rl)\oplus\rl^{10} \\
\hline
17 & (\SU_n,\Upp_n, \Upp_1\cdot\Sp_{n/2})\cdot\SU_2 &
(\cp^n\otimes\cp^2\oplus\rl)\oplus\gt{su}_2 \\
\hline
18 & (\SU_n,\Upp_n, \Upp_1\cdot\Sp_{n/2})\cdot\Upp_2 &
(\cp^n\otimes\cp^2\oplus\rl)\oplus
(\cp^2\oplus\rl)\\
\hline
19 & (\SU_n,\Upp_n, \Upp_1{\cdot}\Sp_{n/2}){\cdot}\SU_2{\cdot}(\SU_n,\Upp_n,
\Upp_1{\cdot}\Sp_{n/2}) & (\cp^n\otimes\cp^2\oplus\rl)\oplus
(\cp^2\otimes\cp^n\oplus\rl)\\
\hline
20 & (\SU_n,\Upp_n, \Upp_1\cdot\Sp_{n/2})\cdot\SU_2\cdot\Upp_4 &
(\cp^n{\otimes}\cp^2{\oplus}\rl){\oplus}
(\cp^2{\otimes}\cp^4{\oplus}\rl){\oplus}\rl^6\\
\hline
21 & \Upp_4\cdot\Upp_2 &
\rl^6\oplus(\cp^4\otimes\cp^2\oplus\rl)\oplus\gt{su}_2\\
\hline
22 & \Upp_4\cdot\Upp_2\cdot\Upp_4 &
\rl^6{\oplus}(\cp^4{\otimes}\cp^2{\oplus}\rl){\oplus}
(\cp^2{\otimes}\cp^4{\oplus}\rl){\oplus}\rl^6\\
\hline
23 & \Upp_1\cdot\Upp_1\cdot\SU_4 &
(\cp^4\oplus\rl)\oplus(\cp^4\oplus\rl)\oplus\rl^6\\
\hline
24 & (\Upp_1\cdot)\SU_4(\cdot\SO_2) & (\cp^4\oplus\rl)\oplus\rl^6\otimes\rl^2\\
\hline
\end{tabular}
\end{center}

\begin{rmk} Suppose $G$ is a complex reductive group. 
Let $X$ be a smooth affine algebraic spherical $G$-variety. 
 By Luna's slice theorem, see \cite{Luna}, 
  we have $X=G\times_H W$, where $H\subset G$ is a reductive subgroup and
  $W$ is a finite-dimensional $H$-module.
As was proved by Knop and Panyushev (private communications),
$H$ is a spherical subgroup of $G$ and 
$W$ is a multiplicity free representation of 
$H_*(\gt g/\gt h)$. Let $K\subset H$ be a maximal 
compact subgroup, $V=(\gt g/\gt h)_\rl$ a $K$-invariant
real form of $\gt g/\gt h$ and $\gt n=W\oplus\rl$ 
a Heisenberg algebra corresponding to $W$. According to 
Theorem~5 and \cite{be-ra},  
the action $K:\gt n\oplus V$ is commutative.  
\end{rmk}

\par\medskip

\begin{center}
{\large \bf 7. Conclusion}
\end{center} 
\nopagebreak

\begin{thm}
Any maximal indecomposable
principal $\Sp_1$-saturated commutative homogeneous space belongs to the
one of the following four classes:

1) affine spherical homogeneous spaces of reductive real Lie groups;

2) spaces corresponding to the rows of Table 2b;

3) homogeneous space $((\rl^n\leftthreetimes\SO_n)\times\SO_n)/\SO_n$,
$((H_n\leftthreetimes\U_n)\times\SU_n)/\Upp_n$,
where the normal subgroups $\SO_n$ and $\SU_n$ of $K$
are diagonally embedded into $\SO_n\times\SO_n$ and
$\SU_n\times\SU_n$, respectively;

4) commutative  homogeneous spaces of Heisenberg type. 
\end{thm}

\begin{proof} Let $X=G/K$ be a commutative homogeneous space. If $G$ is
reductive,
$X$ belongs to the first class.
If $L=K$ then it is  a space of Heisenberg type.

Assume that $G$ is not reductive and $L\ne K$.
Suppose
a simple factor $K_1$ of $K$ has non-trivial projections
onto both $P$ and $L^\diamond$. Then due to condition (1) of the definition of
$\Sp_1$-saturated commutative spaces, $K_1\ne\SU_2$. By theorem 4, $X$ belongs
to the
3-d class. If all simple factors of $K$ are contained in either $P$ or
$L^\diamond$,
then, because $X$ is principal, $P^0/(P^0\cap K)$ is a factor of $X$.
But $X$ is indecomposable
and $G$ is not reductive, so $P^0$ is trivial. Thus, $X$ satisfies condition
($\ast$).

If there is a simple factor $L_1$ of $L$ such that $\pi_1(L_*)\ne K$
and $L_1\subsetneqq K$,
then, according to Theorem 3, $X$ is contained in the second class.
If there is no such factor, then also by Theorem 3, 
$(L, K)$ is a product
of pairs of the type $(\SU_2\times\SU_2\times\SU_2, \SU_2)$,
$(\SU_2\times\SU_2, \SU_2)$ or $(\SU_2, \U_1)$ and a pair
$(K^1, K^1)$, where $K^1$ is a compact Lie group. But
these pairs (except $(K^1, K^1)$) are not allowed in $Sp_1$-saturated
commutative space.
The  second condition
of the definition of $\Sp_1$-saturated commutative space contradicts
the conditions of Lemma 8.
Thus, $L$ would be equal $K$, but this is not the case.
\end{proof}


\begin{thebibliography}{33}
\small

\bibitem{av}
{\sc D.N.~Akhiezer, E.B.~Vinberg,} Weakly symmetric spaces and
spherical varieties, {\it Transformation Groups} {\bf 4} (1999),
3--24.

\bibitem{ap}
{\sc  D.N.~Akhiezer,  D.I.~Panyushev}, Multiplicities in
the branching rules and the complexity
of homogeneous spaces, {\it Moscow Math. J.} {\bf 2} (2002), 17--33.

\bibitem{be-ra}
{\sc C.~Benson, J.~Jenkins, and G.~Ratcliff,} On Gelfand pairs
associated with solvable Lie groups, {\it Trans. AMS} {\bf 321}
(1990),  85--116.

\bibitem{be-ra2}
{\sc C.~Benson and G.~Ratcliff,}
A classification of multiplicity free actions,
{\it Journal of Algebra}, {\bf 181} (1996), no. 1,
152--186.

\bibitem{bg}
{\sc F.A.~Berezin, I.M.~Gelfand, M.I.~Graev, M.A.~Naymark,}
Representations of groups, {\it Uspekhi Mat. Nauk} {\bf 11}, 
no. 6 (1956), 13--40 (in Russian).

\bibitem{br}
{\sc M.~Brion,} Classification des espaces homog\`enes
sph\'eriques, {\it Compositio Math.} {\bf 63} (1987), 189--208.

\bibitem{br2}
{\sc M.~Brion,} Repr\'esentations exceptionnelles des groupes semi-simples,
{\it Ann. Sci. \'Ecole Norm. Sup. (4)},
{\bf 18} (1985), no. 2, 345--387. 

\bibitem{al}
{\sc А.Г.~Элашвили}, {\it Канонический вид и стационарные подалгебры 
точек общего по\-ло\-же\-ния для простых линейных групп Ли}, 
Функц. анализ и его прилож. т.{\bf 6}, вып.~1 (1972), С.~51--62. 
English translation: 
{\sc A.G.~\'Elashvili,}
Canonical form and stationary subalgebras of points in general
position for simple linear Lie groups, {\it Funct. Anal. Appl.}
{\bf 6}, no. 1 (1972), 44--53. 

\bibitem{h}
{\sc S.~Helgason,} Groups and Geometric Analysis, Acad. Press, London, 1984.

\bibitem{Kac}
{\sc Kac, V. G.}, Some remarks on nilpotent orbits,
{\it J. Algebra}, {\bf 64}(1980), no. 1, 190--213.

\bibitem{kr}
{\sc M.~Kr\"amer,} Sph\"arische Untergruppen in kompakten
zusammenh\"angenden Lie\-grup\-pen,  {\it Compositio Math.}
{\bf 38}(1979), 129--153.

\bibitem{kn}
{\sc F.~Knop,} Some remarks on multiplicity free spaces,
 Representation theories and algebraic geometry (Montreal, PQ,
1997), 301--317, NATO Adv. Sci. Inst. Ser. C
 Math. Phys. Sci., 514, Kluwer Acad. Publ., Dordrecht, 1998.

\bibitem{la}
{\sc J.~Lauret,} Gelfand pairs attached to representations 
of compact Lie groups, {\it Transform. Groups} {\bf 5}(2000), 
no. 4, 307--324. 

\bibitem{L}
{\sc A.S~Leahy,} A classification of multiplicity free representations,
{\it J. Lie Theory} {\bf 8}(1998), 367--391.

\bibitem{Luna}
{\sc D.~ Luna,} Slices \'etales,
{\it Bull. Soc. Math. France, Paris, M\'emoire} {\bf 33}(1973),
81--105.

\bibitem{m}
{\sc И.В.~Микитюк,} {\it Об интегрируемости инвариантных гамильтоновых 
систем с одно\-род\-ны\-ми конфигурационными пространствами}, Матем. 
Сб. 1986. Т. 129. С. 514--546. English translation: 
{\sc I.V.~Mikityuk (Mykytyuk),} On the integrability of invariant
Hamiltonian systems with homogeneous configuration spaces, {\it
Math. USSR Sbornik} {\bf 57} (1987), 527--546.  

\bibitem{nish}
{\sc N.~Nishihara,} 
 A geometric criterion for Gelfand pairs associated with
              nilpotent Lie groups,
{\it J. Funct. Anal.}, {\bf 183}\,(2001), no.1, 
148--163.

\bibitem{On}
{\sc A.L.~Onishchik,} Inclusion relations between transitive
transformation groups, {\it Trudy Mosk. Matem. Obshch.} {\bf 11}
(1962), 199--242 (in Russian).

\bibitem{On2}
{\sc A.L.~Onishchik,} Topology of transitive
transformation groups, Johann Ambrosius Barth Verlag GmbH,
Leipzig, 1994. 

\bibitem{dp2}
{\sc D.I.~Panyushev}, Some amazing properties of spherical
nilpotent orbits, {\it Math. Zeitschrift} {\bf 245}(2003), 557--580.

\bibitem{R}
{\sc Л.Г.~Рыбников}, {\it О коммутативности 
слабо коммутативных римановых однородных пространств}, 
Функц. анализ и его прилож. т.{\bf 37}, вып.~2 (2003), С.~41--51. 
English translation: 
{\sc L.G.~Rybnikov}, On commutativity of weakly commutative
Riemannian homogeneous spaces,
{\it Funct. Anal. Appl.} {\bf 37}, no.2
(2003), 114--122.

\bibitem{th}
{\sc E.G.F.~Thomas}, An infinitesimal characterization of Gelfand
pairs, {\it Contemporary Math.} {\bf 26} (1984), 379--385.

\bibitem{Vin}
{\sc Э.Б.~Винберг}, {\it Коммутативные однородные пространства и 
коизотропные действия,} УМН, т.{\bf 56}\,(2001), вып. 1, С.~3--62. 
English translation: 
{\sc E.B.~Vinberg}, Commutative homogeneous
spaces and coisotropic actions, {\it Russian Math. Surveys} 
{\bf 56}, no.1 (2001), 1--60.

\bibitem{Vin2}
{\sc Э.Б.~Винберг}, {\it Коммутативные однородные пространства 
гейзенбергова типа}, Тpуды ММО, т.{\bf 64}\,(2003), С.~54--89. 
English translation:
{\sc E.B.~Vinberg}, Commutative homogeneous
spaces of Heisenberg type, {\it Trans. Moscow Math. Soc.} 
(2003), 47--80.

\bibitem{y2}
{\sc О.С.~Якимова}, {\it О слабо симметрических пространствах 
полупростых групп Ли}, 
Вестник МГУ, Сер. 1, Матем., мех. (2002) no.~2, 57--60.  English translation:
{\sc O.S.~Yakimova},  Weakly symmetric spaces of semisimple Lie
groups, {\it Moscow Univ. Math. Bull.} {\bf 57} (2002), no.\,2, 37-40.

\bibitem{ya3}
{\sc О.С.~Якимова}, {\it О слабо коммутативных однородных 
пространствах}, 
Успехи матем. наук, т.{\bf 57}, вып. 3 (2002), 171--172. 
English translation:
{\sc O.S.~Yakimova}, On weakly commutative homogeneous spaces, {\it
Russian Math. Surveys} {\bf 57}, no.3 (2002), 615--616.

\bibitem{Acta}
{\sc O.S.~Yakimova}, Saturated commutative spaces of Heisenberg type, {\it
Acta Applicandae Mathematica}, {\bf 81}(2004), no.~1, 339--345.
\end{thebibliography}
\end{document}